\documentclass[11pt,reqno,a4paper]{amsart}

\usepackage{lmodern,fbb}
\usepackage[T1]{fontenc}
\usepackage[utf8]{inputenc}
\usepackage{amsmath,amssymb,amsthm,mathtools}

\usepackage[top=2.5cm,bottom=2.4cm,left=2.8cm,right=2.8cm,headsep=0.2in,footskip=13pt]{geometry}

\usepackage{fancyhdr}
\usepackage[shortlabels]{enumitem}
\usepackage[expansion=false]{microtype}
\usepackage{xcolor}
\definecolor{linkteal}{rgb}{0.0,0.42,0.42}
\usepackage[colorlinks=true,linkcolor=linkteal,citecolor=linkteal,urlcolor=linkteal]{hyperref}
\allowdisplaybreaks
\newtheorem{theorem}{Theorem}[section]
\newtheorem{lemma}[theorem]{Lemma}
\newtheorem{proposition}[theorem]{Proposition}
\newtheorem{corollary}[theorem]{Corollary}
\newtheorem{maintheorem}{Theorem}

\theoremstyle{remark}

\newcommand{\R}{\mathbb R}
\newcommand{\Z}{\mathbb Z}

\newcommand{\one}{\mathbf 1}
\newcommand{\BV}{\mathrm{BV}}
\newcommand{\Per}{\mathcal P}
\newcommand{\E}{\mathcal E}
\newcommand{\Div}{\operatorname{div}}
\newcommand{\supp}{\operatorname{supp}}
\newcommand{\diam}{\operatorname{diam}}
\newcommand{\conv}{\operatorname{conv}}
\newcommand{\clip}{\operatorname{clip}_{[0,1]}}

\newcommand{\Lcal}{\mathcal L}
\newcommand{\Fcal}{\mathcal F}
\newcommand{\g}{\mathfrak g}
\newcommand{\aexp}{\alpha}
\title{Sharp isoperimetry in groups of polynomial growth}
\author{Mayukh Mukherjee}

\begin{document}

\begin{abstract}
We study weighted edge boundary and outer vertex boundary on Cayley graphs of polynomial growth of degree $Q>1$. For every finite symmetric generating set, we prove that each exact isoperimetric profile satisfies $I(n)=c n^{(Q-1)/Q}+O(n^{(Q-1)/Q-\beta})$, where $\beta=1/(8Q^2)$. We identify the leading constant with the isoperimetric constant of a horizontal perimeter on the limiting Carnot group, after accounting for the counting density. Its integrand is a norm with polyhedral unit ball, determined by a convex minimization with one real variable for each coset of a torsion-free nilpotent subgroup. For every asymptotically minimizing sequence, we translate its vertex sets so that the characteristic functions of the associated rescaled regions are precompact in $L^1$. The limiting regions minimize continuum perimeter at fixed volume. The exact-cardinality estimate requires no regularity assumption on these regions. For virtually abelian groups, we prove a quantitative stability estimate using a discrete approximation to the Euclidean Wulff body. The relative symmetric difference from a translate is bounded by the square root of the sum of the relative boundary deficit and a cardinality error. We also characterize when bounded additive changes of Cayley metric preserve the leading vertex constant. Failure occurs already in rank two; we construct pairs with the same measured asymptotic cone whose leading constants have an unbounded ratio as the rank grows.
\end{abstract}

\maketitle
\thispagestyle{fancy}
\enlargethispage{5pt}

\section{Introduction}

We study boundary minimization for finite vertex sets in Cayley graphs and its limit under rescaling. Let $\Gamma$ be an infinite group with identity $e$, and let $S=S^{-1}$ be a finite generating set with $e\notin S$. We join each vertex $\gamma$ to $\gamma s$ for $s\in S$, and assign weight $a_s=a_{s^{-1}}>0$ to this edge. For a finite set $A\subset\Gamma$, let $|A|$ denote its cardinality and $\one_A$ its characteristic function. We define its weighted edge boundary $B(A)$ and outer vertex boundary $V(A)$ by
\begin{equation}\label{eq:boundaries}
 B(A)=\sum_{\gamma\in A}\sum_{s\in S}a_s\one_{\Gamma\setminus A}(\gamma s),\qquad V(A)=|AS\setminus A|.
\end{equation}
Thus $B(A)$ counts edges leaving $A$, with their weights, whereas $V(A)$ counts vertices outside $A$ adjacent to $A$. The exact isoperimetric profiles are $I_e(n)=\inf_{|A|=n}B(A)$ and $I_v(n)=\inf_{|A|=n}V(A)$. Multiple copies of a generator in the edge problem are included by adding their weights. We write $\partial=B$ or $V$, and $I=I_e$ or $I_v$ accordingly.

For an abelian group, we regard the generators as vectors in $(\Gamma/\operatorname{Tor}\Gamma)\otimes_{\Z}\R$, where $\operatorname{Tor}\Gamma$ is the torsion subgroup. The limiting problem is anisotropic Euclidean perimeter minimization. Its minimizers are translates and dilates of a convex body: a zonotope, that is, a sum of line segments, for weighted edge boundary; the convex hull of the generator images for outer vertex boundary. These are instances of the Euclidean Wulff construction, which arose from the study of equilibrium crystal shapes. Taylor proved the minimizing and uniqueness statements for bounded sets of finite perimeter; see Fonseca and M\"uller \cite[p.~126]{FonsecaMuller} for the general finite-perimeter formulation and an analytic proof of uniqueness.

Barber and Erde \cite[Theorems~1 and~2]{BarberErde} prove sharp asymptotics for integer lattices; their vertex theorem records an earlier result of Ruzsa. Pansu \cite[Theorem~1]{PansuPeriodic} proves sharp isoperimetric asymptotics for periodic Riemannian metrics on Euclidean space; his Theorem~4 identifies limiting crystal shapes for attractive short-range periodic atomic bond networks in three-dimensional space. Variational limits of lattice perimeter and finite formulas for periodic surface energies are developed in \cite{ChambolleGiacominiLussardi,DelNinPetrache} and \cite[Proposition~2.6 and Theorem~2.7]{ChambolleKreutz}.

Suppose now that $\Gamma$ has polynomial growth of degree $Q>1$: the cardinality of a word ball of radius $R\ge1$ is bounded above and below by positive multiples of $R^Q$. Gromov's theorem gives a nilpotent subgroup of finite index \cite{Gromov}. Its asymptotic cone is a Carnot group: a simply connected nilpotent Lie group whose Lie algebra splits into layers generated by the first under brackets. The stratification is obtained from the associated graded Lie algebra of the nilpotent completion; see \cite[Corollary~1.9]{Breuillard} and \cite{BreuillardLeDonne}. Continuum limits of nonlocal perimeters on Carnot groups are studied in \cite{CarbottiDonPallaraPinamonti}.

The limiting metric and volume measure do not determine which perimeter corresponds to the discrete boundary. We determine it below, separately for edge boundary and outer vertex boundary.

We fix a torsion-free nilpotent normal subgroup $\Lambda\lhd\Gamma$ of finite index $M$, and choose coset representatives $q_1,\ldots,q_M$. Every $\gamma\in\Gamma$ has a unique expression $\gamma=\lambda q_\alpha$ with $\lambda\in\Lambda$ and $1\le\alpha\le M$. Let $N$ be the Malcev completion of $\Lambda$, the connected simply connected nilpotent Lie group in which $\Lambda$ is a lattice. We write $\g$ for the Lie algebra of $N$. Its lower central series determines a graded Lie algebra $V_1\oplus\cdots\oplus V_r$ and its Carnot group $G$, where $r$ is the nilpotency class of $N$. The first layer is the quotient $V_1=\g/[\g,\g]$, and $\pi_1:\g\to V_1$ is the canonical quotient map. We call $V_1$ the horizontal space. We choose a linear splitting to identify $\g$ with the underlying vector space of the graded algebra. The group $N$ need not itself be stratified. Its group operation approaches that of $G$ under rescaling; we describe this convergence in Section~\ref{sec:geometry}.

We normalize Haar measure on $N$ so that $\Lambda$ has covolume one, and transport it to a Haar measure $\mu$ on $G$ through the chosen exponential coordinates. Both measures are the same constant multiple of Lebesgue measure. At scale $h$, we give each vertex mass $h^Q$. The lattice Riemann-sum limit in each coset then gives total counting density $M\mu$. Thus $|A|$ counts vertices in $\Gamma$, whereas $\mu(E)$ measures volume in $G$.

We identify $\lambda q_\alpha$ with $(\lambda,\alpha)$. For each $\alpha$ and $s\in S$, we write uniquely $q_\alpha s=t_{\alpha s}q_{\beta(\alpha,s)}$, with $t_{\alpha s}\in\Lambda$. The generator edge and its horizontal displacement are
\[
 (\lambda,\alpha)\longmapsto(\lambda t_{\alpha s},\beta(\alpha,s)),\qquad
 v_{\alpha s}=\pi_1\log t_{\alpha s}\in V_1.
\]
Here $\log$ is taken in $N$. These edges may join different cosets. For edge boundary, we choose one orientation of each unoriented $\Lambda$-orbit of edges, and denote the finite collection by $\mathcal E_0$. We represent $e\in\mathcal E_0$ by $(1,\alpha(e))\mapsto(t_e,\beta(e))$, where $1$ is the identity of $\Lambda$, and put $v_e=\pi_1\log t_e$. Its weight is denoted by $a_e$.

We now describe the relative boundary levels in the different cosets. We fix a horizontal covector $p\in V_1^*$, that is, a real linear functional on $V_1$. Thus $p(v_{\alpha s})$ and $p(v_e)$ are real numbers obtained by evaluating $p$ on the corresponding horizontal displacements. For $\psi=(\psi_1,\ldots,\psi_M)\in\R^M$, we define the level function
\[
 H_\psi(\lambda,\alpha):=p(\pi_1\log\lambda)+\psi_\alpha.
\]
First-layer coordinates are additive under multiplication in $N$, so the level increases along this edge by $p(v_{\alpha s})+\psi_{\beta(\alpha,s)}-\psi_\alpha$.
For nonzero $p$, the condition $H_\psi\geq0$ defines a discrete set with a separate boundary level in each coset. Changing $\psi$ can change this set. Under rescaling, however, the fixed offsets $\psi_\alpha$ become $h\psi_\alpha$ and tend to zero, so all these sets have the same limiting half-space. Their relative levels can nevertheless affect the leading discrete boundary count. The formulas below choose the least boundary density obtained by varying these levels. The first counts boundary edges with their weights; the second counts exterior boundary vertices, each only once:
\begin{align}
 \phi_e(p)&=\min_{\psi\in\R^M}\sum_{e\in\mathcal E_0}a_e|p(v_e)+\psi_{\beta(e)}-\psi_{\alpha(e)}|,\label{eq:edge-cell}\\
 \phi_v(p)&=\min_{\psi\in\R^M}\sum_{\alpha=1}^M\max\left(0,\max_{s\in S}\{p(v_{\alpha s})+\psi_{\beta(\alpha,s)}-\psi_\alpha\}\right).\label{eq:vertex-cell}
\end{align}
For periodic edge energies, a closely related finite-dimensional minimization and the resulting polyhedrality appear in Pansu \cite[Definition~4 and Proposition~42]{PansuPeriodic}. Here the quotient also records the finite extension of the nilpotent group. The expressions being minimized depend on $\psi$ only through its differences across edges. Adding the same constant to every $\psi_\alpha$ leaves these expressions unchanged, so one may fix the adjustment in a reference coset to be zero.

In the edge formula, the absolute value measures the range of levels for which the endpoints lie on opposite sides of the boundary. Each unoriented edge orbit is included only once. For the vertex formula, we regard $H_\psi\geq0$ as the inside. A vertex of negative level belongs to the exterior boundary if at least one generator step reaches a nonnegative level. For each generator, the possible exterior levels form an interval ending at zero. Their union therefore has length equal to the largest positive level increment, or zero if all increments are nonpositive. This explains the maximum in \eqref{eq:vertex-cell}, rather than a sum which would count some exterior vertices more than once.

The formulas \eqref{eq:edge-cell} and \eqref{eq:vertex-cell} sum contributions from the cosets rather than averaging them. Lemma~\ref{lem:cell-duality} gives a dual description of these formulas and proves that $\phi_e$ and $\phi_v$ are crystalline norms, meaning that their unit balls are polytopes. Propositions~\ref{prop:liminf} and~\ref{prop:recovery} identify them as the limiting perimeter densities: the first proves the lower bound for arbitrary convergent discrete sets, and the second constructs sets attaining the bound in the limit.

For either norm $\phi$, the horizontal perimeter $\Per_\phi(E)$ measures the boundary variation of $E\subset G$ in the directions of $V_1$, weighted by $\phi$; see Section~\ref{sec:perimeter} for the definition. We define the isoperimetric exponent $\aexp$ and the constants $C_\phi,c_\phi$ by
\begin{equation}\label{eq:constant}
 C_\phi=\inf_{0<\mu(E)<\infty}\frac{\Per_\phi(E)}{\mu(E)^{\aexp}},\qquad c_\phi=M^{-\aexp}C_\phi,\qquad \aexp=\frac{Q-1}{Q}.
\end{equation}
The factor $M^{-\aexp}$ converts continuum volume to rescaled cardinality. The cell formulas determine the perimeter density. Finding the sharp constant then means solving an isoperimetric problem on $G$; the proof does not require a classification of its minimizers or regularity of their boundaries.

The isoperimetric inequality on $G$ also has a function version.
The same constant $C_\phi$ is the best constant in the horizontal
$\BV$--Sobolev inequality
\begin{equation}\label{eq:sharp-sobolev}
 C_\phi\|u\|_{Q/(Q-1)}\le\Per_\phi(u)
 \qquad(u\in\BV_H(G)),
\end{equation}
where $\Per_\phi(u)$ is the horizontal variation defined in
Section~\ref{sec:perimeter}. To see this, we apply the isoperimetric
inequality to the level sets of $|u|$. The coarea formula and
Minkowski's inequality give
\[
 \Per_\phi(u)\ge
 C_\phi\int_0^\infty\mu(\{|u|>t\})^{\aexp}\,dt
 \ge C_\phi\|u\|_{Q/(Q-1)}.
\]
\textls[-10]{The constant is sharp because characteristic functions are included
among the functions in} $\BV_H(G)$.

Our first theorem shows that the isoperimetric problem on $G$
determines the leading term of the discrete profile. It also shows
that near-minimizers on $\Gamma$ approach isoperimetric regions in $G$
after translation and rescaling.

\begin{maintheorem}\label{thm:main}
Suppose $\Gamma$ has polynomial growth of degree $Q>1$. For each of
the boundaries in \eqref{eq:boundaries}, let $\phi$ be the corresponding
norm in \eqref{eq:edge-cell} or \eqref{eq:vertex-cell}. With
$\beta=1/(8Q^2)$, there is a constant $C<\infty$, depending on the
group and the generating data, such that
\begin{equation}\label{eq:profile-rate}
 |I(n)-c_\phi n^{\aexp}|\le C n^{\aexp-\beta}\qquad(n\ge1).
\end{equation}
Moreover, let $A_n\subset\Gamma$ satisfy $|A_n|=n$ and
$\partial(A_n)/n^{\aexp}\to c_\phi$. We associate to $A_n$ one rescaled
region in $G$ for each coset, as in Section~\ref{sec:geometry}.
After suitable left translations of $A_n$ by elements of $\Lambda$,
the characteristic functions of the corresponding regions are
precompact in $L^1$. Every limit is $\one_E$ in each coset, where
$E\subset G$ is an isoperimetric region with $M\mu(E)=1$ and
$\Per_\phi(E)=c_\phi$.
\end{maintheorem}

A sharper exponent, depending on the nilpotency class $r$, is given in \eqref{eq:explicit-beta}. At scale $h=n^{-1/Q}$, the $M$ associated regions have total volume one and approach the same region $E$. The leading constant and this convergence statement are independent of the subgroup, representatives, splitting, and bounded fundamental domain.

The distinction between $I(n)$ and its greatest nondecreasing minorant $\inf_{m\ge n}I(m)$ is the subject of a question of Gromov \cite[p.~501]{GromovEntropy}. A comparison up to fixed multiplicative constants is proved for elementary amenable groups, among other classes, in \cite[Corollary~6.12]{MukherjeeExactMonotone}. For groups of polynomial growth, Corollary~\ref{cor:monotone-folner} below gives the sharper ratio $1+O(n^{-\beta})$, as well as the leading constant and a power error for the corresponding F\o lner function.

Theorem~\ref{thm:main} concerns asymptotically optimal sets. The next theorem applies to arbitrary finite sets $A\subset\Gamma$, separately for $h^{Q-1}B(A)$, the rescaled weight of boundary edges, and $h^{Q-1}V(A)$, the rescaled number of outside vertices adjacent to $A$. Their limiting perimeters are determined by $\phi_e$ and $\phi_v$, respectively.

Here $\Gamma$-convergence means that every convergent sequence satisfies the lower bound and that, for every limiting set, there are approximating discrete sets attaining this bound. Convergence is strong in $L^1(G)$ in each coset. Tightness means uniformly small total $L^1$ mass outside sufficiently large compact subsets of $G$.

\begin{maintheorem}\label{thm:gamma}
Under nilpotent rescaling, each functional $h^{Q-1}\partial$ $\Gamma$-converges in strong $L^1(G)^M$ to
\[
 (\one_E,\ldots,\one_E)\longmapsto\Per_\phi(E),
\]
with value $+\infty$ on all other tuples. Families of finite sets with bounded rescaled volume and boundary are locally precompact, and they are precompact in strong $L^1$ when tight. For every measurable $E\subset G$ of finite volume and finite $\phi$-perimeter, and every prescribed sequence $h\downarrow0$, there are finite sets $A_h\subset\Gamma$ whose rescaled parts converge to $E$ in every coset and whose rescaled boundaries satisfy
\[
 h^{Q-1}\partial(A_h)\longrightarrow\Per_\phi(E).
\]
\end{maintheorem}

For virtually abelian groups, a small boundary deficit gives a quantitative estimate at finite cardinality. Here $A$ is a finite subset of $\Gamma$, and the comparison sets are discrete approximations to minimizers of anisotropic perimeter in Euclidean space. Suppose $\Lambda$ is free abelian of rank $d\ge2$. We identify $G$ with $\R^d$ and retain the Haar normalization and coset representatives above. Let $W\subset\R^d$ be the Wulff body of $\phi$, defined by $\phi(p)=\max_{x\in W}p(x)$. We set
\[
 W_*=(M\mu(W))^{-1/d}W,\qquad
 \mathcal W_n=\{\lambda q_\alpha:\lambda\in\Lambda\cap n^{1/d}W_*,\
 1\le\alpha\le M\}.
\]
Thus $M\mu(W_*)=1$, and $\mathcal W_n\subset\Gamma$ is a lattice approximation to the volume-normalized continuum minimizer $W_*\subset G$. It satisfies $|\mathcal W_n|=n+O(n^{(d-1)/d})$. We write $t_+=\max(t,0)$ for the positive part of a real number $t$.

\begin{maintheorem}\label{thm:abelian-stability}
For either weighted edge boundary or outer vertex boundary there is $C<\infty$, depending on the fixed data, such that every nonempty finite $A\subset\Gamma$, with $|A|=n$, satisfies
\begin{equation}\label{eq:abelian-stability}
 \inf_{\kappa\in\Lambda}\frac{|A\triangle\kappa\mathcal W_n|}{n}
 \le C\left(\delta_\phi(A)+n^{-1/(2d)}\right)^{1/2},\qquad
 \delta_\phi(A)=\left(\frac{\partial(A)}{c_\phi n^{(d-1)/d}}-1\right)_+.
\end{equation}
\end{maintheorem}

The left-hand side of \eqref{eq:abelian-stability} is the proportion of vertices on which $A$ differs from a translated Wulff approximation. Unlike the limiting assertion in Theorem~\ref{thm:main}, this estimate applies to each finite set. The comparison set describes the shape and need not itself have asymptotically optimal boundary: the relative levels in the cosets may still need adjustment.

For integer lattices, Barber--Erde--Keevash--Roberts \cite[Theorems~1.1 and~1.2]{BarberErdeKeevashRoberts} prove the corresponding quantitative vertex estimate and qualitative edge stability. Theorem~\ref{thm:abelian-stability} gives a quantitative estimate for general positive edge weights and extends both estimates to finite extensions. We use the quantitative Wulff inequality of Figalli--Maggi--Pratelli \cite[Theorem~1.1]{FigalliMaggiPratelli}; for its application to discrete energies on periodic lattices, see Cicalese--Leonardi \cite{CicaleseLeonardi}.

For a torsion-free nilpotent group we may take $M=1$. The formulas reduce to
\begin{equation}\label{eq:one-coset}
 \phi_e(p)=\sum_j a_j|p(X_j)|,\qquad \phi_v(p)=\max_{s\in S}p(X_s),\qquad X_s=\pi_1\log s,
\end{equation}
where $j$ indexes one chosen generator from each inverse pair, with weight $a_j$ and horizontal image $X_j$. Consequently the leading constants depend on the associated graded group, the covolume, and the horizontal images of the generators. Adding generators with zero horizontal image changes neither leading constant. These conclusions hold for Carnot lattices and for lattices in non-stratified nilpotent groups.

Two word metrics that differ by at most a fixed constant have the same limit under rescaling, but their leading vertex-boundary constants can differ. For virtually abelian groups, Theorem~\ref{thm:holonomy-characterization} determines exactly when this can happen. The constant is unchanged under every such bounded change of metric if and only if the conjugation action of the finite quotient on $\Lambda\otimes_{\Z}\R$ is contained in $\{I,-I\}$. Otherwise, there are generating sets $S\subset T$ such that
\[
 d_T\le d_S\le d_T+C,\qquad c_{v,S}<c_{v,T}.
\]
Corollary~\ref{cor:holonomy-unbounded} gives examples of rank $d$ for which $c_{v,T}/c_{v,S}\ge(d!)^{1/(2d)}$. In particular, this ratio is unbounded as $d$ grows. We also compute the leading constants for three wallpaper groups. For the standard generators of the Heisenberg group, we give explicit upper and lower bounds for the leading edge constant. We construct sets of every cardinality $n$ whose edge boundary is at most $(8/3^{3/4})n^{3/4}+Cn^{1/4}$. We also prove a minimization result for sets between Lipschitz graphs over a square.

We briefly describe the proofs of Theorems~\ref{thm:main} and~\ref{thm:gamma}. To obtain the lower bound in the variational limit, we use duality with divergence-free edge functions on the finite quotient. For the matching upper bound, we start with smooth functions on the limiting group, evaluate them at rescaled lattice points, and adjust their values separately in each coset. Convolution on the left and on the right relates the discrete boundary to horizontal perimeter and gives compactness. For this part, we only need bounded fundamental cells; their boundaries need not be regular. A partition of the graph extends the local lower bound to arbitrary finite sets. Since $0<\aexp<1$, a set split into pieces of positive limiting mass has too large a boundary to be nearly optimal. Finally, subadditivity allows us to construct sets of every prescribed cardinality. The same argument applies to edge and vertex boundaries.

\section{Nilpotent coordinates and the finite quotient}\label{sec:geometry}

The purpose of these coordinates is to make rescaled generator differences converge to horizontal derivatives on $G$. Let $\g$ be the Lie algebra of the completion $N$ fixed above. Its lower central series is $\g^1=\g$, $\g^{k+1}=[\g,\g^k]$; the nilpotency class $r$ is determined by $\g^r\ne0$ and $\g^{r+1}=0$. We choose complements $W_k$ with $\g^k=W_k\oplus\g^{k+1}$. The associated graded algebra is $\g_\infty=\bigoplus_{k=1}^rV_k$, where $V_k=\g^k/\g^{k+1}$, with the bracket induced by that of $\g$. Projection identifies each chosen complement $W_k$ with $V_k$. The splitting therefore identifies $\g=\bigoplus_{k=1}^rW_k$ with the underlying vector space of $\g_\infty$; it need not identify their brackets. Using exponential coordinates of the first kind, we may consequently regard $N$ and $G$ as the same vector space with different group products. For $Z_k\in W_k$, we define the dilation $d_h$ and homogeneous dimension $Q$ by
\[
 d_h\left(\sum_k Z_k\right)=\sum_k h^k Z_k,\qquad Q=\sum_k k\dim W_k.
\]
The number $Q$ is the growth degree introduced above. For $h>0$, we express multiplication in $N$ in these rescaled coordinates by
\begin{equation}\label{eq:transport-product}
 x*_h y=d_h\bigl(d_h^{-1}x\cdot d_h^{-1}y\bigr).
\end{equation}
The products $*_h$ extend polynomially to $h=0$, where $*_0$ is the product of $G$. In particular, they converge smoothly on compact sets. Haar measure is the same constant multiple of Lebesgue measure for all these products. Inversion is $x\mapsto-x$ in these coordinates.

We fix a bounded Borel left fundamental domain $D$ for $\Lambda$ in $N$, with $\mu(D)=1$. We put $z_\lambda=d_h\lambda$ and $D_h=d_hD$, using the exponential-coordinate convention for $\lambda\in\Lambda$. The cells
\[
 C_{\lambda,h}=d_h(\lambda D)=z_\lambda*_h D_h
\]
partition the underlying manifold up to null sets and have measure $h^Q$. We call $z_\lambda$ the lattice point associated with $C_{\lambda,h}$; it need not belong to the cell. For an array of coefficients $f=(f_{\lambda\alpha})$, we define the cell function $u_h^\alpha$ to have value $f_{\lambda\alpha}$ on $C_{\lambda,h}$. We also write $f_h^\alpha$ for this function. For a finite set $A\subset\Gamma$, we take $f_{\lambda\alpha}=\one_A(\lambda q_\alpha)$. In particular,
\begin{equation}\label{eq:mass}
 \sum_\alpha\|u_h^\alpha\|_1=h^Q|A|.
\end{equation}
We fix a norm on $V_1$ and the corresponding Carnot--Carath\'eodory distance, defined by minimizing lengths of horizontal curves. We choose a homogeneous norm $|x|_0$ on $G$ comparable to the Carnot--Carath\'eodory distance from the identity, and write $B_R=\{x:|x|_0<R\}$. For $\gamma\in\Gamma$, let $|\gamma|_S$ be its word length and $d_S(\gamma,\gamma')=|\gamma^{-1}\gamma'|_S$ the word distance. We write $B_S(\gamma,R)$ for the word ball and $\diam_S$ for the corresponding diameter. In estimates, $C,c>0$ denote constants whose values may change between occurrences. They may depend on the fixed group, fundamental domain, generating set, and weights, but are independent of $h$ unless stated otherwise.

The next lemma makes the convergence of the rescaled products quantitative. Its main assertion is that a fixed lattice step, after rescaling, acts on a smooth function like a small horizontal displacement in $G$.

\begin{lemma}\label{lem:coordinates}
For $v\in V_1$, let $X_v$ be the left-invariant horizontal vector field on $G$ with value $v$ at the identity.
For $R\ge1$, $0<h\le1$, and fixed $t\in\Lambda$, we put $v_t=\pi_1\log t$. The constants below may depend on $t$, and are uniform when $t$ ranges over a fixed finite subset of $\Lambda$. All coordinate derivatives and $C^2$ norms are taken in the fixed exponential coordinates. If $b$ is $C^2$ on a neighborhood of $B_{CR}$, then
\begin{equation}\label{eq:increment}
 b(x*_h d_ht)=b(x)+hX_{v_t}b(x)+O\left(h^2R^{2r}\|b\|_{C^2(B_{CR})}\right)\qquad(|x|_0\le R).
\end{equation}
If $z_\lambda\in B_R$, then $\sup_{y\in C_{\lambda,h}}|y-z_\lambda|_{\rm Eucl}\le Ch(1+R)^{r-1}$; in particular these cells have diameter at most this bound after increasing $C$. Products of points in $B_R$ lie in $B_{CR}$. The first and second coordinate derivatives of the map $(x,y)\mapsto x*_h y$ on $B_R\times B_R$ are bounded by $CR^{r-1}$ and $CR^{\max(r-2,0)}$, respectively. Moreover,
\begin{equation}\label{eq:word-growth}
 C^{-1}|\lambda|_0-C\le|\lambda|_S\le C|\lambda|_0+C,\qquad cR^Q\le|B_S(e,R)|\le CR^Q\quad(R\ge1),
\end{equation}
where $|\lambda|_0$ means the homogeneous norm of $\log\lambda$ in the chosen splitting. In particular the Cayley graph is doubling: every ball of radius $2R$ is covered by a uniformly bounded number of balls of radius $R$.
\end{lemma}

\begin{proof}
We give a coordinate in $W_k$ weight $k$, and denote the coordinate weights by $w_i\le r$. For a multi-index $I$, we write $|I|_w=\sum_iw_iI_i$. Each monomial of the transported Baker--Campbell--Hausdorff product has a real coefficient $c$ and the form $ch^{w_i-|I|_w-|J|_w}x^Iy^J$, with $|I|_w+|J|_w\le w_i$. This gives the derivative bounds. Products stay in $B_{CR}$ by
\[
 d_{1/R}(x*_h y)=(d_{1/R}x)*_{h/R}(d_{1/R}y).
\]
We substitute $y=d_ht$. A term of first order in $h$ can involve only the first-layer component $v_t$ of $\log t$, and its coefficient comes from the graded product $*_0$. Differentiating $b(x*_0(\tau v_t))$ at $\tau=0$ gives $X_{v_t}b(x)$, since right multiplication by this small step produces the left-invariant field. All remaining increments contain at least $h^2$. We increase $C$ so that the Euclidean segment between $x$ and $x*_h d_ht$ lies in $B_{CR}$; the coordinate bounds above allow this even though a homogeneous ball need not be Euclidean convex. Taylor's theorem on this segment gives \eqref{eq:increment}. Substitution of a bounded element of $D$ gives the cell estimate.

The word-length comparison in \eqref{eq:word-growth} is the usual homogeneous quasi-norm comparison for nilpotent lattices; it remains valid on a finite-index subgroup with the ambient word metric. Guivarc'h's comparison is recalled in \cite[Theorem~2.7]{Breuillard}. The volume bounds and their extension to virtually nilpotent groups follow from \cite[Corollary~2.9 and Remark~2.10]{Breuillard}, or by counting the fundamental cells. Finite index does not change the exponent.
\end{proof}

We extend the boundary counts from sets to functions so that we can approximate and smooth them. The coarea identity will relate these extensions to the boundaries of their level sets. For a finitely supported real array $f=(f_{\lambda\alpha})$, with $\lambda\in\Lambda$ and $1\le\alpha\le M$, we regard $f$ as a function on $\Gamma$ by $f(\lambda q_\alpha)=f_{\lambda\alpha}$. We define the convex extensions of the two boundaries by
\begin{align}
 \Lcal_e(f)&=\sum_{\lambda\in\Lambda}\sum_{e\in\mathcal E_0}a_e|f_{\lambda t_e,\beta(e)}-f_{\lambda,\alpha(e)}|,\label{eq:convex-edge}\\
 \Lcal_v(f)&=\sum_{\lambda,\alpha}\max\left(0,\max_{s\in S}\{f_{\lambda t_{\alpha s},\beta(\alpha,s)}-f_{\lambda\alpha}\}\right).\label{eq:convex-vertex}
\end{align}
For characteristic functions these equal $B(A)$ and $V(A)$, respectively. For Euclidean lattices, these convex extensions and their coarea identities were developed by Chambolle--Giacomini--Lussardi \cite{ChambolleGiacominiLussardi}. We write $\Lcal$ for the corresponding convex energy. Both energies satisfy the exact coarea identity
\begin{equation}\label{eq:discrete-coarea}
 \Lcal(f)=\int_0^1\partial(\{f>t\})\,dt\qquad(0\le f\le1).
\end{equation}
Indeed, in the vertex case the integrand at a vertex is one precisely when its value is at most $t$ and some neighbor's value is greater than $t$. In either case, we write $\E_h(f)=h^{Q-1}\Lcal(f)$ for the rescaled energy. For nonnegative $f$, its rescaled mass is $m_h(f)=h^Q\sum_{\lambda,\alpha}f_{\lambda\alpha}$.

Differences between the values at two vertices can be estimated by adding differences along a path of generator edges. The following estimates will in particular force the different coset components to have the same limit.

\begin{lemma}\label{lem:paths}
For every finitely supported real-valued $f$, either energy controls the sum of the absolute coefficient differences across all directed generator edges, up to a fixed constant. In particular,
\begin{align}
 h^Q\sum_\lambda|f_{\lambda\alpha}-f_{\lambda\beta}|&\le Ch\E_h(f),\label{eq:coset-sync}\\
 h^Q\sum_\lambda|f_{\lambda\kappa,\alpha}-f_{\lambda\alpha}|&\le Ch(1+|\kappa|_S)\E_h(f)\qquad(\kappa\in\Lambda).\label{eq:right-shift}
\end{align}
In particular, \eqref{eq:coset-sync} says that, when the rescaled energies stay bounded, $\|u_h^\alpha-u_h^\beta\|_1\to0$ as $h\downarrow0$.
\end{lemma}

\begin{proof}
For the edge energy this follows from the positive weights. For the vertex energy, we sum the positive part of every directed increment. Each is at most the maximum at its initial vertex, and reversing an edge accounts for its negative part. Thus the total absolute generator variation is bounded by $2|S|\Lcal_v(f)$.

We join $\lambda q_\alpha$ to $\lambda q_\beta$ using a fixed word for $q_\alpha^{-1}q_\beta$. For each position in the word, the map from $\lambda$ to the initial vertex of the edge is injective. Telescoping and summing proves \eqref{eq:coset-sync}. For \eqref{eq:right-shift}, we use a word for $q_\alpha^{-1}\kappa q_\alpha$, whose length is at most $|\kappa|_S+C$.
\end{proof}

\section{The effective horizontal perimeter}\label{sec:perimeter}

The cell formula gives a cost for each horizontal covector, rather than the perimeter of a whole set. To identify this cost by duality, we combine edge increments so that every relative level $\psi_\alpha$ cancels. On the finite quotient graph, the sum of coefficients leaving each vertex must then equal the sum entering it. These balanced coefficients give the lower bound in Proposition~\ref{prop:liminf}; the minimizing levels give the matching construction in Proposition~\ref{prop:recovery}.

The quotient graph has vertices $1,\ldots,M$. We retain distinct edges
with the same endpoints, since their displacements may differ. In the
edge problem we use the oriented edges in $\mathcal E_0$; in the vertex
problem we use all directed generator edges. For a real-valued
function $F$ on the relevant edge collection, we define
\[
 (\Div_0F)_\alpha
 =\sum_{\alpha(e)=\alpha}F_e
  -\sum_{\beta(e)=\alpha}F_e.
\]
Thus $\Div_0F=0$ is exactly the cancellation condition. Indeed,
\[
 \sum_e F_e\bigl(\psi_{\beta(e)}-\psi_{\alpha(e)}\bigr)
 =-\sum_\alpha\psi_\alpha(\Div_0F)_\alpha.
\]

The coefficient bounds come directly from the two cell formulas.
An absolute value allows a signed coefficient with $|F_e|\le a_e$.
A maximum at a vertex allows nonnegative coefficients whose sum is
at most one. Accordingly, we put
\begin{align*}
 \Fcal_e
 &=\{F:|F_e|\le a_e,\ \Div_0F=0\},\\
 \Fcal_v
 &=\{F:F_{\alpha s}\ge0,\
             \sum_sF_{\alpha s}\le1,\ \Div_0F=0\}.
\end{align*}
These are compact convex polytopes. In either case, we write $\Fcal$ for
the corresponding polytope and define
\[
 T(F)=\sum_eF_ev_e,\qquad W=T(\Fcal)\subset V_1.
\]
We can choose the minimizing relative levels continuously as the normal varies. The polyhedral form of the cell cost lets us choose them linearly on finitely many cones. We fix Euclidean norms on $V_1^*$ and $\R^M$, and let $L(p,\psi)$ denote the expression minimized in the corresponding cell formula.

\begin{lemma}\label{lem:cell-duality}
The minima in \eqref{eq:edge-cell}--\eqref{eq:vertex-cell} are attained, and
\begin{equation}\label{eq:dual-cell}
 \phi(p)=\max_{F\in\Fcal}p(T(F))=\max_{w\in W}p(w).
\end{equation}
The set $W$ is a symmetric full-dimensional polytope. In particular, $\phi$ is a norm on $V_1^*$.

There is a positively homogeneous, piecewise linear map $\psi:V_1^*\to\R^M$ such that
\begin{equation}\label{eq:corrector-lipschitz}
 \psi_1(p)=0,\qquad L(p,\psi(p))=\phi(p),\qquad
 |\psi(p)-\psi(q)|\le C|p-q|.
\end{equation}
In particular, $\psi(0)=0$ and $|\psi(p)|\le C|p|$.
\end{lemma}

\begin{proof}
Adding the same constant to all components of $\psi$ leaves $L$
unchanged. We remove this freedom by setting $\psi_1=0$.
Since the quotient graph is connected, every component of $\psi$ can
then be estimated by adding differences along a path from vertex $1$.
On each edge,
\[
 |\psi_{\beta(e)}-\psi_{\alpha(e)}|
 \le |p(v_e)+\psi_{\beta(e)}-\psi_{\alpha(e)}|+|p(v_e)|.
\]
In the edge case the positive weights control the first term by $L$.
In the vertex case, reversing a directed edge changes the sign of
its increment. Summing positive parts therefore controls the sum of
absolute increments by $2|S|L$. It follows in both cases that
\[
 |\psi|\le C\bigl(L(p,\psi)+|p|\bigr)\qquad(\psi_1=0).
\]
Thus a normalized minimizing sequence is bounded, and continuity
gives a minimizer. Since $L(p,0)\le C|p|$, every normalized minimizer
satisfies $|\psi|\le C|p|$.

To obtain the dual formula, we use the elementary identities
\[
 a|z|=\max_{|F|\le a}Fz,\qquad
 \max(0,z_1,\ldots,z_\ell)
 =\max_{\substack{F_j\ge0\\\sum_jF_j\le1}}\sum_jF_jz_j.
\]
These turn each cell formula into a finite linear program.
In its dual, the coefficient of $\psi_\alpha$ is
$-(\Div_0F)_\alpha$. Minimization over the shifts
therefore imposes $\Div_0F=0$. The remaining expression is
$\sum_eF_ep(v_e)=p(T(F))$. Linear programming duality proves
\eqref{eq:dual-cell}.

The image $W$ is a polytope because $T$ is linear. In the edge case,
$F\mapsto-F$ shows that $W$ is symmetric. For the vertex case, let
$\bar e$ be the reverse of a directed generator edge $e$. The symmetry
of $S$ gives $v_{\bar e}=-v_e$. We put $F'_e=F_{\bar e}$.
The new outgoing sum at a vertex is the old incoming sum,
which equals the old outgoing sum because $\Div_0F=0$. The new
coefficients are therefore admissible, and their displacement is
$-T(F)$. Hence $W$ is symmetric in this case as well.

To see that $W$ has full dimension, we take $\lambda\in\Lambda$ and
follow a generator path from $q_1$ to $\lambda q_1$. Its projection
to the quotient is a closed walk. The sum of its horizontal
displacements is $\pi_1\log\lambda$, since projection to the
abelianization turns multiplication in $N$ into addition.
We count its signed edge traversals in the edge problem, and its
directed traversals in the vertex problem. At every vertex the
number entering equals the number leaving. A sufficiently small
positive multiple of these counts therefore belongs to $\Fcal$.
Thus $W$ contains a positive multiple of $\pi_1\log\lambda$.
These vectors span $V_1$ as $\lambda$ ranges over the lattice.
Since $W$ is convex and symmetric, it contains a neighborhood of
the origin. Its support function \eqref{eq:dual-cell} is consequently
a norm.

To choose the minimizers, we triangulate the boundary of the polytope $\{p:\phi(p)\le1\}$ so that each simplex lies in one facet and the triangulations agree on shared faces. At each vertex $p_i$, we choose a minimizer $\psi_i$ with $(\psi_i)_1=0$. On the cone over a simplex, we write uniquely $p=\sum_i a_i p_i$, with $a_i\ge0$, and set $\psi(p)=\sum_i a_i\psi_i$. The function $\phi$ is linear on this cone, so $\phi(p)=\sum_i a_i$. Joint convexity and positive homogeneity of $L$ give
\[
 L(p,\psi(p))\le\sum_i a_iL(p_i,\psi_i)
 =\sum_i a_i=\phi(p).
\]
The reverse inequality follows from the definition of $\phi$. The formulas agree on common faces and therefore define a continuous piecewise linear map. Since there are finitely many linear pieces, the map is globally Lipschitz. This proves \eqref{eq:corrector-lipschitz}.
\end{proof}

For a smooth function, we apply the cell cost to its horizontal differential at each point and integrate over the group. For the indicator of a set of finite horizontal perimeter, the horizontal derivative is a measure. Weighting this measure by the same cell cost defines the perimeter of the set.

We fix a basis $(e_i)$ of $V_1$ and its dual basis $(e^i)$, and use the corresponding Euclidean norm. For a smooth function $u$, its horizontal differential is the covector $\nabla_Hu=\sum_i(X_{e_i}u)e^i$. For $u\in L^1(G)$, we take these derivatives in the distributional sense relative to $\mu$. We write $u\in\BV_H(G)$ when $D_Hu=\sum_i(X_{e_i}u)e^i$ is a finite $V_1^*$-valued measure, and denote its total variation by $|D_Hu|$. We define
\[
 \Per_\phi(u)=\int_G\phi\left(\frac{dD_Hu}{d|D_Hu|}\right)\,d|D_Hu|,\qquad \Per_\phi(E)=\Per_\phi(\one_E).
\]
Equivalently, by \eqref{eq:dual-cell},
\begin{equation}\label{eq:perimeter-dual}
 \Per_\phi(u)=\sup\left\{-\int_Gu\sum_eX_{v_e}F_e\,d\mu:F\in C_c^\infty(G;\Fcal)\right\}.
\end{equation}
To see that this gives the full anisotropic variation, we choose preimages in $\Fcal$ of the vertices of $W$. Smooth convex combinations of these preimages approximate measurable choices against the finite measure $|D_Hu|$. The resulting fields are admissible and attain the usual dual supremum in the limit.

We use lower semicontinuity, coarea, and smooth approximation for horizontal variation; the differential formulation and approximation theorem are treated in \cite[Theorems~2.4 and~3.1]{AmbrosioGhezziMagnani}. For the anisotropic perimeter, smooth approximation also follows directly from convolution. We choose a nonnegative smooth compactly supported kernel $\rho$ with integral one, and put $\rho_\varepsilon(x)=\varepsilon^{-Q}\rho(d_{1/\varepsilon}x)$. Unless a subscript indicates otherwise, convolution uses the graded product on $G$, with $(\rho*u)(x)=\int_G\rho(y)u(y^{-1}x)\,d\mu(y)$. Left convolution commutes with each horizontal derivative, so
\begin{equation}\label{eq:convolution-perimeter}
 \Per_\phi(\rho_\varepsilon*u)\le\Per_\phi(u).
\end{equation}
Cutoff at large radii and lower semicontinuity therefore give smooth compactly supported approximants with convergent $\phi$-variation. For $0\le u\le1$, the horizontal coarea formula is
\[
 \Per_\phi(u)=\int_0^1\Per_\phi(\{u>t\})\,dt.
\]

\section{Quantitative smoothing}\label{sec:smoothing}

We turn a function on the discrete group into a smooth function on $G$. The smooth function must remain close in $L^1$, and its perimeter must not exceed the discrete energy apart from a small error. The coefficient of the discrete energy must be exactly one; a fixed multiplicative loss would lose the sharp limiting constant.

The path estimate controls smoothing by right translations. The sharp perimeter estimate is more naturally expressed using left convolution. We first use right smoothing to control the approximation, then pass to left smoothing. For the energy estimate, we smooth the test fields. Only boundedness of the fundamental cells is needed, not a bound on their perimeter. Throughout this section, $\rho$ is nonnegative, smooth, compactly supported in a fixed homogeneous ball, invariant under inversion, and has integral one. We write $\clip(t)=\min(1,\max(0,t))$.

\begin{lemma}\label{lem:smoothing}
Let $\theta=1/(r+1)$, and let $\delta_z$ denote the unit point mass at $z\in G$. For each fixed $R\ge1$, there is $C_R<\infty$ with the following property. Suppose $0\le f\le1$ is finitely supported and its associated lattice points lie in $\{|x|_0\le R\}$. Let $m=m_h(f)$, $E=\E_h(f)$, and
\[
 \nu_h=\frac{h^Q}{M}\sum_{\lambda,\alpha}f_{\lambda\alpha}\delta_{z_\lambda},\qquad v=\rho_\varepsilon*\nu_h.
\]
For $0<h,\varepsilon<1$,
\begin{align}
 \|\clip(v)-u_h^\alpha\|_1&\le C_R\bigl(\varepsilon^\theta+h\varepsilon^{-2r}\bigr)(E+m),\label{eq:smooth-approx}\\
 \Per_\phi(v)&\le E+C_R h\varepsilon^{-2r}(E+m).\label{eq:smooth-sharp}
\end{align}
In \eqref{eq:smooth-approx}, $u_h^\alpha$ is the cell function of $f$ in the indicated coset.
\end{lemma}

\begin{proof}
Throughout this proof, constants may depend on $R$. We first prove the $L^1$ estimate. For Euclidean coordinate derivatives, the kernel satisfies $\|D^j\rho_\sigma\|_\infty\le C_j\sigma^{-Q-jr}$ and $\|D^j\rho_\sigma\|_1\le C_j\sigma^{-jr}$. By Lemma~\ref{lem:coordinates}, the translated kernels $K(x,y)=\rho_\sigma(x*_0y^{-1})$ and $K(x,y)=\rho_\sigma(y^{-1}*_h x)$ have first and second derivative bounds in $L^1$ at most $C\sigma^{-r}$ and $C\sigma^{-2r}$ on the relevant supports. For finitely supported real coefficients $b_\lambda$ with $b_\lambda=0$ unless $|z_\lambda|_0\le R$, Lemma~\ref{lem:coordinates} gives $\sup_{y\in C_{\lambda,h}}|y-z_\lambda|_{\rm Eucl}\le Ch$. We apply the mean-value theorem between the associated point and each point of its cell, using the derivative bound on a fixed enlargement of the ball. This gives
\begin{equation}\label{eq:kernel-cell}
 \left\|h^Q\sum_\lambda b_\lambda K(\,\cdot\,,z_\lambda)-\sum_\lambda b_\lambda\int_{C_{\lambda,h}}K(\,\cdot\,,y)\,d\mu(y)\right\|_1
 \le Ch\sigma^{-r}h^Q\sum_\lambda|b_\lambda|.
\end{equation}
The same comparison gives Riemann sums for smooth functions.

Right convolution for the transported product lets us apply \eqref{eq:right-shift} to estimate the $L^1$ error:
\[
 w_\sigma(x)=h^Q\sum_\lambda f_{\lambda\alpha}\rho_\sigma(z_\lambda^{-1}*_h x),\qquad
 w_\sigma(z_\eta)=h^Q\sum_{\kappa\in\Lambda}f_{\eta\kappa,\alpha}\rho_\sigma(d_h\kappa^{-1}).
\]
The sum $Z=h^Q\sum_\kappa\rho_\sigma(d_h\kappa^{-1})$ is close to one. Indeed, for $h\le\sigma$, the cells meeting the support of $\rho_\sigma$, or with their associated points in that support, have union contained in $B_{C\sigma}$. This follows by dilation by $(\sigma+h)^{-1}$, for which the transported product has parameter $h/(\sigma+h)\le1$. The union has volume at most $C\sigma^Q$. Its point-to-cell distances are at most $Ch$, so the bound $\|D\rho_\sigma\|_\infty\le C\sigma^{-Q-r}$ gives
\[
 \sum_\lambda\int_{C_{\lambda,h}}
 |\rho_\sigma(y)-\rho_\sigma(z_\lambda)|\,d\mu(y)
 \le Ch\sigma^{-r},\qquad |Z-1|\le Ch\sigma^{-r}.
\]
Here inversion reindexes the sum defining $Z$.
When this error is at most $1/2$, we normalize the nonnegative coefficients by $Z$. Every $\kappa$ in their support satisfies $|\kappa|_S\le C(\sigma/h+1)$, so \eqref{eq:right-shift} controls the difference between the normalized lattice values and $f_{\eta\alpha}$ by $C(\sigma+h)E$. The normalization contributes at most $|Z-1|m$.

To compare $w_\sigma$ with its lattice values on the cells, we translate each kernel center to the identity using $*_h$. This exactly reindexes both the cells and their associated points:
\[
 z_\kappa^{-1}*_h C_{\eta,h}=C_{\kappa^{-1}\eta,h},\qquad
 z_\kappa^{-1}*_h z_\eta=z_{\kappa^{-1}\eta}.
\]
The preceding integrated kernel estimate therefore bounds this error by $Ch\sigma^{-r}m$. If $h>\sigma$ or the normalization error bound exceeds $1/2$, the direct estimate $\|w_\sigma-u_h^\alpha\|_1\le2m$ gives the same conclusion after increasing $C$. Differentiating the convolution kernel gives the second estimate below. Thus
\begin{equation}\label{eq:opposite-smooth}
 \begin{split}
 \|w_\sigma-u_h^\alpha\|_1&\le C(\sigma+h)E+Ch\sigma^{-r}m,\\
 \|\nabla_{\rm Eucl}w_\sigma\|_1&\le C\sigma^{-r}m.
 \end{split}
\end{equation}

To pass to left convolution, we take $a\in\supp\rho_\varepsilon$ and interpolate by $a_t=\exp_G(t\log_G a)$, $0\le t\le1$. Each map $x\mapsto a_t^{-1}*_0x$ preserves Haar measure. Its coordinate velocity on the relevant support is bounded by $C\varepsilon$. Consequently the fundamental theorem of calculus gives
\[
 \|w_\sigma(a^{-1}*_0\,\cdot\,)-w_\sigma\|_1
 \le C\varepsilon\|\nabla_{\rm Eucl}w_\sigma\|_1.
\]
We insert $w_\sigma$ between each translate of $u_h^\alpha$ and $u_h^\alpha$, average over $a$, and use \eqref{eq:opposite-smooth}. We obtain
\[
 \|\rho_\varepsilon*u_h^\alpha-u_h^\alpha\|_1
 \le C\bigl(\sigma+h+h\sigma^{-r}+\varepsilon\sigma^{-r}\bigr)(E+m).
\]
The choice $\sigma=\varepsilon^{1/(r+1)}$ balances the terms $\sigma$ and $\varepsilon\sigma^{-r}$. The error in replacing the averaged cell density by $\nu_h$ is at most $Ch\varepsilon^{-r}m$, by \eqref{eq:kernel-cell}; replacing it by $u_h^\alpha$ introduces an error at most $ChE$, by \eqref{eq:coset-sync}. These estimates prove \eqref{eq:smooth-approx} before truncation to $[0,1]$, since the remaining powers of $\varepsilon^{-1}$ are at most $2r$. Truncation cannot increase the distance from a function valued in $[0,1]$.

We now prove the sharp energy estimate by testing against the dual fields in \eqref{eq:perimeter-dual}. We fix $F\in C_c^\infty(G;\Fcal)$ and set
\[
 \widetilde F_e(z)=\int_G\rho_\varepsilon(x*_0z^{-1})F_e(x)\,d\mu(x).
\]
This remains in $\Fcal$ pointwise: its constraints are linear and convex, and the kernel has integral one. Its $C^2$ norm on the relevant compact region is bounded by $C\varepsilon^{-2r}$, independently of derivatives of $F$. With $a=x*_0z^{-1}$, right invariance of Haar measure gives
\[
 \widetilde F_e(z)=\int_G\rho_\varepsilon(a)F_e(a*_0z)\,d\mu(a).
\]
We differentiate along $z*_0\exp_G(tv_e)$. The corresponding curve in the integrand is $(a*_0z)*_0\exp_G(tv_e)$, so its derivative is $X_{v_e}F_e(a*_0z)$, with no change of sign. Thus
\begin{equation}\label{eq:test-convolution}
 X_{v_e}\widetilde F_e(z)=\int_G\rho_\varepsilon(x*_0z^{-1})X_{v_e}F_e(x)\,d\mu(x).
\end{equation}
Discrete summation by parts gives
\begin{align*}
 J_h&:=h^{Q-1}\sum_{\lambda,e}\widetilde F_e(z_\lambda)(f_{\lambda t_e,\beta(e)}-f_{\lambda,\alpha(e)})\\
 &=h^{Q-1}\sum_{\lambda,\alpha}f_{\lambda\alpha}\left(\sum_{\beta(e)=\alpha}\widetilde F_e(z_\lambda*_h d_ht_e^{-1})-\sum_{\alpha(e)=\alpha}\widetilde F_e(z_\lambda)\right).
\end{align*}
In the vertex case this sum is over all directed edges. Admissibility gives $J_h\le E$ for either energy. The zeroth-order terms cancel because $\Div_0\widetilde F=0$. The step $t_e^{-1}$ has horizontal part $-v_e$, so \eqref{eq:increment} first gives
\[
 J_h=-h^Q\sum_{\lambda,e}f_{\lambda,\beta(e)}X_{v_e}\widetilde F_e(z_\lambda)
 +O\bigl(h\varepsilon^{-2r}m\bigr).
\]
We now replace each $f_{\lambda,\beta(e)}$ by $\bar f_\lambda=M^{-1}\sum_\alpha f_{\lambda\alpha}$. The resulting error is bounded by \eqref{eq:coset-sync} times the first-derivative bound on $\widetilde F$. Hence
\[
 J_h=-h^Q\sum_{\lambda,e}\bar f_\lambda X_{v_e}\widetilde F_e(z_\lambda)+O\bigl(h\varepsilon^{-2r}m+h\varepsilon^{-r}E\bigr).
\]
By \eqref{eq:test-convolution}, the leading term is $-\int v\sum_eX_{v_e}F_e\,d\mu$. Taking the supremum in \eqref{eq:perimeter-dual} proves \eqref{eq:smooth-sharp}.
\end{proof}

\begin{corollary}\label{cor:compactness}
As $h\downarrow0$, families of characteristic functions of the associated regions in $G$ with uniformly bounded rescaled mass and energy are locally precompact in $L^1(G)^M$. Tight families are precompact in strong $L^1(G)^M$, and every limit has the form $(\one_E,\ldots,\one_E)$.
\end{corollary}

\begin{proof}
We first suppose that the associated lattice points lie in a fixed ball. For fixed $\varepsilon>0$, the convolutions $v_h=\rho_\varepsilon*\nu_h$ from Lemma~\ref{lem:smoothing} have a common compact support and uniformly bounded first derivatives. They are therefore precompact in $L^1$, as are their truncations $\clip(v_h)$, since truncation is $L^1$-continuous. By \eqref{eq:smooth-approx},
\[
 \limsup_{h\downarrow0}\max_\alpha\|u_h^\alpha-\clip(v_h)\|_1
 \le C\varepsilon^\theta.
\]
We let $\varepsilon\downarrow0$ and take a diagonal subsequence. This argument applies to all coefficients in $[0,1]$.

For a compact set $K$, we choose $\chi\in C_c^\infty(G)$ with $0\le\chi\le1$ and $\chi=1$ on a neighborhood of $K$, and put $g_{\lambda\alpha}=\chi(z_\lambda)f_{\lambda\alpha}$. The increment estimate gives, for both energies,
\[
 m_h(g)\le m_h(f),\qquad \E_h(g)\le\E_h(f)+C_\chi m_h(f).
\]
Indeed, we expand each product difference, bound the original increment using $0\le\chi\le1$, and reindex the term containing the value of $f$ at the terminal vertex. Since $C_{\lambda,h}=z_\lambda*_h d_hD$ and $D$ is bounded, the cell functions of $g$ agree with those of $f$ on $K$ for all sufficiently small $h$. The preceding argument gives compactness on $K$, and a diagonal subsequence over a compact exhaustion gives local compactness.

By tightness, local $L^1$ convergence implies global $L^1$ convergence. Local strong limits of indicators are indicators, and \eqref{eq:coset-sync} makes the coset limits identical.
\end{proof}

\section{The variational limit}\label{sec:gamma}

We now show that the perimeter defined in Section~\ref{sec:perimeter} is the limiting discrete boundary. First, no sequence of discrete sets approaching $E$ can have a smaller limiting boundary. Then we construct discrete sets approaching $E$ whose boundaries attain this value: we smooth the characteristic function, adjust its lattice values in each coset, and take a common level set.

\begin{proposition}\label{prop:liminf}
Let $A_h\subset\Gamma$ be finite, with cell functions $u_h^\alpha$. If $u_h^\alpha\to u^\alpha$ in $L^1(G)$ for all $\alpha$, then
\[
 \liminf_{h\downarrow0}\E_h(\one_{A_h})\ge\Per_\phi(E)
\]
when all $u^\alpha=\one_E$; otherwise the liminf is $+\infty$.
\end{proposition}

\begin{proof}
If the liminf is infinite there is nothing to prove. Otherwise we pass to a subsequence realizing the finite liminf; its energies are bounded. The path estimate makes the limits equal, and strong convergence makes the common limit an indicator. For an arbitrary $F\in C_c^\infty(G;\Fcal)$, discrete duality gives $\E_h(\one_{A_h})\ge J_h$, where summation by parts yields
\begin{align*}
 J_h&=h^{Q-1}\sum_{\lambda,\alpha}f_{\lambda\alpha}
 \left(\sum_{\beta(e)=\alpha}F_e(z_\lambda*_h d_ht_e^{-1})
       -\sum_{\alpha(e)=\alpha}F_e(z_\lambda)\right)\\
 &=-h^Q\sum_{\lambda,e}f_{\lambda,\beta(e)}X_{v_e}F_e(z_\lambda)+o(1).
\end{align*}
Here $f=\one_{A_h}$. The zeroth-order terms cancel by $\Div_0F=0$, and the error tends to zero by \eqref{eq:increment} on a fixed compact enlargement of the support of $F$. We replace $f_{\lambda,\beta(e)}$ by the average $\bar f_\lambda=M^{-1}\sum_\alpha f_{\lambda\alpha}$. By \eqref{eq:coset-sync}, the resulting error is $O(h\E_h(f))=o(1)$. The cell mesh estimate and strong convergence now give
\[
 \liminf_h\E_h(\one_{A_h})\ge-\int_G\one_E\sum_eX_{v_e}F_e\,d\mu.
\]
Taking the supremum over $F$ in \eqref{eq:perimeter-dual} proves the assertion. This also covers a limit set of infinite perimeter.
\end{proof}

For the upper bound, we evaluate a smooth function at the lattice points and adjust its values separately in each coset. We use the exact minimizing map $\psi$ from Lemma~\ref{lem:cell-duality}. For $w\in C_c^\infty(G)$, $0\le w\le1$, we set
\begin{equation}\label{eq:corrected-samples}
 f_{\lambda\alpha}=\clip\left(w(z_\lambda)+h\psi_\alpha(\nabla_Hw(z_\lambda))\right).
\end{equation}
The adjustments tend to zero with $h$ but attain the optimal cell cost. The following estimates will also be used for the power remainder.

\begin{lemma}\label{lem:quant-recovery}
Let $K_0$ be a fixed compact set. There are $C,h_0>0$, depending only on $K_0$ and the fixed data, such that, for $0<h\le h_0$, the lattice values \eqref{eq:corrected-samples} associated with $w\in C_c^\infty(K_0)$, $0\le w\le1$, satisfy
\begin{align}
 \E_h(f)&\le\Per_\phi(w)+Ch\|w\|_{C^2},\label{eq:quant-recovery-energy}\\
 \sum_\alpha\|f_h^\alpha-w\|_1&\le Ch\|w\|_{C^1}.\label{eq:quant-recovery-mass}
\end{align}
All derivatives are taken in the fixed exponential coordinates.
\end{lemma}

\begin{proof}
We put $p(x)=\nabla_Hw(x)$ and $g_\alpha(x)=w(x)+h\psi_\alpha(p(x))$. These functions vanish outside $K_0$, since $\psi(0)=0$. We apply \eqref{eq:increment} to the smooth function $w$. For the added term, we use only the Lipschitz estimate in \eqref{eq:corrector-lipschitz}:
\[
 |\psi_\beta(p(x*_h d_ht_e))-\psi_\beta(p(x))|
 \le Ch\|w\|_{C^2}.
\]
Consequently each generator difference of $g$, divided by $h$, equals its increment in the cell formula up to $Ch\|w\|_{C^2}$. Absolute values and finite maxima preserve this error bound. The cell cost is exactly
$L(p(x),\psi(p(x)))=\phi(p(x))$, whose spatial Lipschitz constant is at most $C\|w\|_{C^2}$. The point-to-cell estimate therefore gives
\[
 \left|h^Q\sum_\lambda\phi(p(z_\lambda))-
       \int_G\phi(p(x))\,d\mu(x)\right|
 \le Ch\|w\|_{C^2}.
\]
Only cells in a fixed compact enlargement of $K_0$ contribute, and their total volume is bounded independently of $h$. Summing the generator errors proves the energy estimate, and truncation to $[0,1]$ cannot increase either energy. Finally, $|g_\alpha-w|\le Ch\|w\|_{C^1}$ and the point-to-cell estimate give the mass estimate, also after truncation.
\end{proof}

\begin{proposition}\label{prop:recovery}
For every measurable $E\subset G$ of finite volume and finite $\phi$-perimeter, and every prescribed sequence $h\downarrow0$, there are finite sets $A_h\subset\Gamma$ whose cell functions satisfy $u_h^\alpha\to\one_E$ strongly in $L^1$ for every $\alpha$ and
\[
 \limsup_h\E_h(\one_{A_h})\le\Per_\phi(E).
\]
\end{proposition}

\begin{proof}
We choose $w_j\in C_c^\infty(G)$, $0\le w_j\le1$, with $w_j\to\one_E$ in $L^1$ and $\Per_\phi(w_j)\to\Per_\phi(E)$. For each fixed $j$, Lemma~\ref{lem:quant-recovery} makes the cell functions approach $w_j$ and bounds their energy by $\Per_\phi(w_j)+o_h(1)$. Along the prescribed scales, we choose $j=j(h)\to\infty$ slowly enough that these errors vanish. The resulting cell functions approach $\one_E$ in every coset, with energy at most $\Per_\phi(E)+o(1)$.

We write $f^\alpha=f_h^\alpha$ for these cell functions and $e_h=\sum_\alpha\|f^\alpha-\one_E\|_1$. We choose $\delta_h\downarrow0$ slowly enough that $e_h/\delta_h\to0$. For any $0<\delta<1/2$ and every common level $t\in[\delta,1-\delta]$, one has
\begin{equation}\label{eq:level-distance}
 \sum_\alpha\|\one_{\{f^\alpha>t\}}-\one_E\|_1\le\delta^{-1}\sum_\alpha\|f^\alpha-\one_E\|_1.
\end{equation}
The coarea identity supplies one level $t_h\in[\delta_h,1-\delta_h]$ whose rescaled boundary is at most $\E_h(f)/(1-2\delta_h)$. We use this same level in every coset. Equation~\eqref{eq:level-distance} gives strong convergence, and the displayed boundary estimate gives the required upper bound along the prescribed scales.
\end{proof}

Propositions~\ref{prop:liminf} and~\ref{prop:recovery}, together with Corollary~\ref{cor:compactness}, prove Theorem~\ref{thm:gamma}. For a set of finite volume and infinite perimeter, we approximate in $L^1$ by finite unions of cells. The upper-bound condition is then automatic because its right-hand side is infinite.

\paragraph{Independence of the auxiliary choices.}
For fixed $k_\alpha\in\Lambda$, changing the representatives to $q'_\alpha=k_\alpha q_\alpha$ replaces an edge displacement by
\[
 v'_e=v_e+\pi_1\log k_{\alpha(e)}-\pi_1\log k_{\beta(e)}.
\]
We set $\psi'_\alpha=\psi_\alpha+p(\pi_1\log k_\alpha)$. Then
$p(v'_e)+\psi'_{\beta(e)}-\psi'_{\alpha(e)}
 =p(v_e)+\psi_{\beta(e)}-\psi_{\alpha(e)}$, so both cell costs are unchanged. The corresponding lattice points differ by a fixed lattice step, which vanishes after rescaling on compact sets.

Two choices of splitting induce the identity on every quotient $\g^k/\g^{k+1}$. Their coordinate change therefore has only higher-layer corrections; after conjugation by $d_h$, each such correction has a positive power of $h$. The rescaled changes tend locally uniformly to the canonical graded identification and preserve the transported Haar measure. Thus they describe the same limiting sets and perimeter.

For two bounded fundamental domains attached to the same lattice points, cell integrals against a continuous compactly supported test function differ by $o(1)$, while total masses agree. This weak comparison suffices here: if $0\le u_h\le1$ converges weakly to $\one_E$ and $\int u_h\to\mu(E)$, then
\[
 \|u_h-\one_E\|_1=\int u_h+\mu(E)-2\int_Eu_h\longrightarrow0.
\]
The convergence of the last integral follows by approximating $\one_E$ in $L^1$ by compactly supported continuous functions and using $0\le u_h\le1$. Hence both choices give the same strong limit.

Finally let $\Lambda'\lhd\Gamma$ be contained in $\Lambda$, with $[\Lambda:\Lambda']=d$. Its normalized Haar measure and number of cosets are
\[
 \mu'=\mu/d,\qquad M'=dM.
\]
We choose representatives $r_1,\ldots,r_d$ for $\Lambda/\Lambda'$ and group the $d$ refined cosets $\Lambda' r_jq_\alpha$ over each old coset. Measured using $\mu$, the average $d^{-1}\sum_j u'^{\alpha,j}_h$ has exactly the mass of $u_h^\alpha$. Their integrals against compactly supported continuous functions differ by $o(1)$, as follows by comparing the atomic sums and using the bounded shifts $r_j$. If $u_h^\alpha\to\one_E$ strongly, the preceding argument makes this average converge strongly to $\one_E$. Since every summand lies in $[0,1]$,
\[
 \left\|d^{-1}\sum_j u'^{\alpha,j}_h-\one_E\right\|_1
 =d^{-1}\sum_j\|u'^{\alpha,j}_h-\one_E\|_1,
\]
so each refined cell function converges strongly to $\one_E$. Conversely, strong convergence of the refined cell functions gives convergence of the old ones by the same comparison. This identifies the convergence notions without using the perimeter conclusion.

Both discrete energies are unchanged, so the two variational limits give the same numerical perimeter functional, $\Per'_{\phi}(E)=\Per_\phi(E)$. Consequently
\[
 C'_\phi=d^{\aexp}C_\phi,\qquad
 c'_\phi=(dM)^{-\aexp}C'_\phi=M^{-\aexp}C_\phi=c_\phi.
\]
Passing to a common finite-index normal subgroup compares any two initial choices.
\section{Sharp profiles and concentration}\label{sec:global}

A finite set may be spread over a large part of the graph. We cut it into pieces of controlled diameter while adding little boundary, then apply the local lower bound to each piece. The same partition prevents loss of mass for near-minimizers; subadditivity fills the gaps between the cardinalities supplied by approximation.

\begin{lemma}\label{lem:coarse-partition}
There is $c_0>0$ such that $\partial(A)\ge c_0|A|^{\aexp}$ for every nonempty finite set. For every integer $L\ge2$ and finite $A$, there is a partition $A=\bigsqcup_i A_i$ with
\begin{equation}\label{eq:partition}
 \diam_S A_i\le4L,\qquad \sum_i\partial(A_i)\le\partial(A)+C|A|/L.
\end{equation}
\end{lemma}

\begin{proof}
We choose $R_A\le C|A|^{1/Q}$ with $|B_S(e,R_A)|\ge2|A|$. Since $\sum_g|A\cap Ag|=|A|^2$, some $g$ in this ball satisfies $|A\triangle Ag|\ge|A|$. We write $g=s_1\cdots s_m$ with $m\le R_A$, and put $f=\one_A$. Telescoping gives
\begin{align*}
 |A\triangle Ag|
 &=\sum_x|f(xg)-f(x)|\\
 &\le\sum_{j=1}^m\sum_x
 |f(xs_1\cdots s_j)-f(xs_1\cdots s_{j-1})|\\
 &=\sum_{j=1}^m\sum_y|f(ys_j)-f(y)|\le CR_A\partial(A).
\end{align*}
The last equality reindexes $y=xs_1\cdots s_{j-1}$; it does not require right invariance of the word metric. The generator variation bound is Lemma~\ref{lem:paths}. This proves the coarse inequality.

For the partition, we choose a maximal $L$-separated net, fix an ordering of its points, and choose one common integer radius uniformly from $\{L,\ldots,2L\}$. Maximality ensures that these balls cover the graph. We assign each vertex to the first net point whose ball contains it. Each cluster lies in one ball, so has diameter at most $4L$; it need not be connected.

If neighboring vertices receive different assignments, some relevant net ball contains one and not the other. For each net point this permits at most two choices of the common radius. Only net points within $2L+1$ of the edge matter. Their $L/3$-balls are disjoint and lie in a ball of comparable radius, so the growth bounds in \eqref{eq:word-growth} give a uniform bound on their number. Thus an edge is separated with probability at most $C/L$. Summing over edges incident to $A$ gives a realization with at most $C|A|/L$ separated edges.

We intersect the clusters with $A$. In the edge case each separated internal edge is added twice, with its weight. In the vertex case, an exterior vertex $z\notin A$ is already counted once in $V(A)$. Any additional appearance in a different piece's boundary can be charged to a cut edge from $z$ to that piece, since at most one piece belongs to the cluster of $z$. A vertex $z\in A$ appearing in another piece's exterior boundary is likewise charged to an internal cut edge. Each cut edge receives only a bounded number of charges. This proves \eqref{eq:partition}.
\end{proof}

We apply the coarse inequality to the approximating sets from Proposition~\ref{prop:recovery}. Passing to the limit gives
\[
 \Per_\phi(E)\ge c_0(M\mu(E))^{\aexp},\qquad C_\phi\ge c_0M^{\aexp}>0.
\]
Infinite perimeter causes no difficulty. A bounded smooth set shows $C_\phi<\infty$. Dilations give $\Per_\phi(d_tE)=t^{Q-1}\Per_\phi(E)$ and $\mu(d_tE)=t^Q\mu(E)$.

\begin{lemma}\label{lem:local-sharp}
Let $R<\infty$, $\delta>0$, and $\varepsilon>0$ be fixed. For all sufficiently small $h$, if
\[
 \diam_S A\le R/h,\qquad \delta\le h^Q|A|\le1,
\]
then
\begin{equation}\label{eq:local-sharp}
 h^{Q-1}\partial(A)\ge(c_\phi-\varepsilon)(h^Q|A|)^{\aexp}.
\end{equation}
\end{lemma}

\begin{proof}
If $\varepsilon\ge c_\phi$, the conclusion is immediate. Otherwise we suppose that a sequence violates it. Its energies are bounded by $c_\phi$, and after taking a subsequence its masses converge to $m\in[\delta,1]$. We translate by an element of $\Lambda$ so that one vertex lies in the fixed set of coset representatives. The diameter bound and Lemma~\ref{lem:coordinates} place all cell supports in one fixed ball. Compactness therefore gives strong convergence to a common $\one_E$, with $M\mu(E)=m\ge\delta>0$. Proposition~\ref{prop:liminf} yields
\[
 \liminf h^{Q-1}\partial(A)\ge C_\phi\mu(E)^{\aexp}=c_\phi m^{\aexp},
\]
whereas violation of \eqref{eq:local-sharp} gives a limsup at most $(c_\phi-\varepsilon)m^{\aexp}$. The positive mass makes these bounds incompatible.
\end{proof}

\begin{proposition}\label{prop:sharp-law}
One has $I(n)/n^{\aexp}\to c_\phi$.
\end{proposition}

\begin{proof}
For the lower bound we take $|A|=n$, $h=n^{-1/Q}$, and partition at scale $L=\lfloor R/h\rfloor$. Let $m_i=h^Q|A_i|$, so $\sum_i m_i=1$ and
\[
 h^{Q-1}\partial(A)\ge\sum_i h^{Q-1}\partial(A_i)-C/R+o(1).
\]
We choose $\delta>0$ with $c_0\delta^{\aexp-1}\ge c_\phi$. For $m_i\ge\delta$ we apply Lemma~\ref{lem:local-sharp} with $4R$ in place of $R$; for $m_i<\delta$ the coarse inequality gives $c_0m_i^{\aexp}\ge c_\phi m_i$. Since $m_i^{\aexp}\ge m_i$, summation gives a lower bound $c_\phi-\varepsilon-C/R+o(1)$. We let $R\to\infty$ and $\varepsilon\downarrow0$.

To fill the gaps between the cardinalities supplied by approximation, we use subadditivity: we translate one of two finite sets on the left by an element of $\Lambda$ far enough that no edge joins the sets and their outer boundaries are disjoint. Boundary is unchanged by this translation, so $I(n+m)\le I(n)+I(m)$. We choose $E$ with $M\mu(E)=1$ and $\Per_\phi(E)\le c_\phi+\varepsilon$. Proposition~\ref{prop:recovery}, applied at $h=1/k$, gives finite sets of cardinalities $N_k$ with $N_k/k^Q\to1$ and $I(N_k)\le(c_\phi+\varepsilon+o(1))N_k^{\aexp}$. We set $N_k^*=\max_{j\le k}N_j$. Then $N_k^*/k^Q\to1$, successive distinct values of $N_k^*$ have ratio tending to one, and the same energy estimate holds for each $N_k^*$.

These values first give $I(n)\le Cn^{\aexp}$ for every $n$. For all sufficiently large $n$, the ratio property supplies a value $N_k^*\in[n/2,n]$. We apply subadditivity and repeat on the remainder until it falls below a fixed threshold. The resulting powers sum geometrically because $\aexp>0$; the finitely many possible final remainders are absorbed by increasing $C$. For the sharp estimate we choose the largest value $N_*=N_k^*\le n$. Then $n-N_*=o(n)$, and
\[
 I(n)\le I(N_*)+C(n-N_*)^{\aexp}\le(c_\phi+\varepsilon+o(1))n^{\aexp}.
\]
Letting $\varepsilon\downarrow0$ completes the proof.
\end{proof}

\begin{proposition}\label{prop:concentration}
If $|A_n|=n$ and $\partial(A_n)/n^{\aexp}\to c_\phi$, then, at scale $h=n^{-1/Q}$, left translations by elements of $\Lambda$ make their cell functions tight. For every strong subsequential limit $(\one_E,\ldots,\one_E)$, one has $M\mu(E)=1$ and $\Per_\phi(E)=c_\phi$.
\end{proposition}

\begin{proof}
With $h=n^{-1/Q}$, we show that, for each $\eta>0$, there is $R_\eta<\infty$ such that a word ball of radius $R_\eta/h$ contains at least $(1-\eta)n$ vertices for all sufficiently large $n$. The radius may depend on $\eta$. Otherwise the partition splits the mass into pieces, and strict concavity of $t^{\aexp}$ makes their total boundary too large. We fix $0<\eta<1/2$ and put $q_\eta=(1-\eta)^{\aexp-1}>1$. If no ball of radius $R/h$ contains $(1-\eta)n$ vertices, we partition at the integer scale $L=\lfloor R/(4h)\rfloor$; for small $h$, $L\ge2$, and every piece has diameter at most $R/h$, hence mass $m_i\le1-\eta$. We choose $\delta$ so that $c_0\delta^{\aexp-1}\ge c_\phi q_\eta$. The local sharp inequality for larger pieces and the coarse inequality for smaller ones give
\[
 h^{Q-1}\partial(A_n)\ge(c_\phi-\varepsilon)q_\eta-C/R+o(1).
\]
We choose $\varepsilon$ small and then $R$ large to make this strictly larger than $c_\phi$. This is impossible for infinitely many $n$. Thus, for each $\eta$, a ball of radius $R_\eta/h$ contains $(1-\eta)n$ vertices eventually.

We choose centers for $\eta=1/4$. A ball containing at least $3n/4$ vertices intersects any ball containing at least $(1-\eta)n$ vertices, with $\eta<1/4$. After enlarging the radii, these same centers give concentration for every $\eta$. We replace each center by its $\Lambda$ coordinate, a change of bounded word distance, and translate that coordinate to the identity. The coordinate comparison transfers concentration to tightness of the cells. After passing to a subsequence, Corollary~\ref{cor:compactness} gives strong convergence. The mass identity and the liminf inequality give $M\mu(E)=1$ and $\Per_\phi(E)\le c_\phi$; the reverse inequality is its definition.
\end{proof}

\begin{corollary}\label{cor:bounded-minimizer}
The infimum in \eqref{eq:constant} is attained by a bounded set of any prescribed positive volume.
\end{corollary}

Existence and boundedness hold more generally in sub-Finsler nilpotent groups by \cite[Theorems~1.1 and~4.9]{PozueloExistence}. Here existence follows from discrete approximation, and boundedness from a continuum truncation argument.

\begin{proof}
We choose approximate discrete minimizers and apply Propositions~\ref{prop:sharp-law}--\ref{prop:concentration}. This gives a continuum minimizer. We dilate it so that $\mu(E)=1$, hence $\Per_\phi(E)=C_\phi$.

We choose a positive degree-one homogeneous function $\varrho$ that is smooth away from the identity and comparable to $|\cdot|_0$. Its horizontal gradient is homogeneous of degree zero and bounded. We put $t(R)=\mu(E\setminus\{\varrho<R\})$ and write $d\mu=c_\mu\,dx$ in the chosen coordinates. On compact annuli, the Euclidean coarea formula makes $t$ locally absolutely continuous. Indeed, the weighted Euler identity implies $\nabla_{\rm Eucl}\varrho\ne0$ away from the identity, and for almost every $R$,
\[
 -t'(R)=c_\mu\int_{E\cap\{\varrho=R\}}\frac{1}{|\nabla_{\rm Eucl}\varrho|}\,d\mathcal H^{\dim G-1}.
\]
Cutting $E$ along this level adds two copies of the cut boundary. Its anisotropic horizontal density is
$c_\mu\phi(\nabla_H\varrho)/|\nabla_{\rm Eucl}\varrho|$, so the added perimeter is at most $C(-t'(R))$. The isoperimetric inequality on the two pieces therefore gives
\[
 C_\phi\bigl((1-t(R))^{\aexp}+t(R)^{\aexp}\bigr)
 \le\Per_\phi(E\cap\{\varrho<R\})+\Per_\phi(E\setminus\{\varrho<R\})
 \le C_\phi+C(-t'(R))
\]
for almost every $R$. For small positive $t$, the left side minus $C_\phi$ is at least $(C_\phi/2)t^{\aexp}$. Thus $-t'\ge ct^{\aexp}$ wherever the tail is small and positive. If it stayed positive for all large $R$, integration of
$(t^{1-\aexp})'\le-c(1-\aexp)$ would make it negative at a finite radius. Hence $t$ vanishes at a finite radius. Replacing $E$ on a null set gives a bounded minimizer; dilation gives every positive volume.
\end{proof}

\section{A power remainder}\label{sec:rate}

For the lower bound, we balance smoothing error against the diameter of the partition pieces. In Lemma~\ref{lem:smoothing}, balancing $\varepsilon^\theta$ and $h\varepsilon^{-2r}$, where $\theta=1/(r+1)$, gives an error of order $h^\gamma$, with
\[
 \gamma=\frac{\theta}{2r+\theta}=\frac1{2r(r+1)+1}.
\]

\begin{lemma}\label{lem:quant-local}
For each fixed $\delta>0$ there is $C_\delta<\infty$ such that, for small $h$ and every $R\ge1$,
\begin{equation}\label{eq:quant-local}
 \begin{gathered}
 \diam_S A\le R/h,\qquad \delta\le m=h^Q|A|\le1\\
 \Longrightarrow\quad
 h^{Q-1}\partial(A)\ge\bigl(c_\phi-C_\delta R^{\aexp}h^{\aexp\gamma}\bigr)m^{\aexp}.
 \end{gathered}
\end{equation}
\end{lemma}

\begin{proof}
We put $E_h=h^{Q-1}\partial(A)$; we may assume $E_h\le2c_\phi$. We translate as in Lemma~\ref{lem:local-sharp}, and apply Lemma~\ref{lem:smoothing} at the smaller scale $h'=h/R$. The associated lattice points then lie in one fixed ball, independently of $R$. The energy and mass at this scale are $R^{1-Q}E_h$ and $R^{-Q}m$.

Let $w'$ be the truncated smoothing at scale $h'$. We put $w(x)=w'(d_{1/R}x)$, and write $u=u_h^1$. The cells satisfy $C_{\lambda,h'}=d_{1/R}C_{\lambda,h}$ exactly. Volume scales by $R^Q$, and horizontal perimeter by $R^{Q-1}$. Lemma~\ref{lem:smoothing} therefore gives
\begin{align*}
 \|w-u\|_1
 &\le C\bigl(\varepsilon^\theta+(h/R)\varepsilon^{-2r}\bigr)(RE_h+m),\\
 \Per_\phi(w)
 &\le E_h+C(h/R)\varepsilon^{-2r}(E_h+m/R).
\end{align*}
Truncation does not increase perimeter. We take $\varepsilon=(h/R)^{1/(2r+\theta)}$. Since $E_h$ and $m$ are bounded,
\[
 \|w-u\|_1\le Ch^\gamma R^{1-\gamma},\qquad
 \Per_\phi(w)\le E_h+C(h/R)^\gamma.
\]
The path estimate gives $|M\|u\|_1-m|\le ChE_h$. We put $p=1/\aexp$. Since $u$ is an indicator and $0\le w\le1$, one has $\|w-u\|_p\le\|w-u\|_1^{\aexp}$. We apply the sharp Sobolev inequality \eqref{eq:sharp-sobolev} to obtain
\begin{align*}
 E_h
 &\ge c_\phi m^{\aexp}
   -C\bigl(h^{\aexp}+h^{\aexp\gamma}R^{\aexp(1-\gamma)}+(h/R)^\gamma\bigr)\\
 &\ge c_\phi m^{\aexp}-CR^{\aexp}h^{\aexp\gamma}.
\end{align*}
Finally we use $m^{\aexp}\ge\delta^{\aexp}$ to obtain the stated coefficient.
\end{proof}

\begin{proposition}\label{prop:rate-lower}
With $\beta=\aexp\gamma/[Q(1+\aexp)]$, one has
$I(n)\ge c_\phi n^{\aexp}-Cn^{\aexp-\beta}$.
\end{proposition}

\begin{proof}
We take $|A|=n$ and choose $\delta$ with $c_0\delta^{\aexp-1}\ge c_\phi$. We set $h=n^{-1/Q}$ and $R=h^{-\aexp\gamma/(1+\aexp)}$. For small $h$, the coefficient in Lemma~\ref{lem:quant-local} is positive. We apply the partition from Proposition~\ref{prop:sharp-law} and this lemma to the larger pieces, with $4R$ in place of $R$. The smaller pieces contribute at least $c_\phi$ times their mass. Thus
\[
 h^{Q-1}\partial(A)\ge c_\phi-C\bigl(R^{\aexp}h^{\aexp\gamma}+R^{-1}\bigr).
\]
Both errors are $O(h^{\aexp\gamma/(1+\aexp)})$, which proves the assertion.
\end{proof}

For the upper bound, we smooth a bounded continuum minimizer and apply Lemma~\ref{lem:quant-recovery}. We first shrink the minimizer slightly so that a common level set has at most $n$ vertices. We add the missing vertices far away by subadditivity. All supports before this final addition stay in a fixed compact set.

\begin{lemma}\label{lem:continuum-smooth-error}
Let $\rho$ be the nonnegative smooth compactly supported kernel of integral one fixed in Section~\ref{sec:smoothing}. If $u\in\BV_H(G)$ is supported in a fixed ball and $0<\varepsilon\le1$, then, with $\theta=1/(r+1)$,
\begin{equation}\label{eq:continuum-smooth-error}
 \|\rho_\varepsilon*u-u\|_1\le C\varepsilon^{\theta}\bigl(|D_Hu|(G)+\|u\|_1\bigr).
\end{equation}
For bounded $u$, it also satisfies $\|\rho_\varepsilon*u\|_{C^2}\le C\varepsilon^{-2r}\|u\|_\infty$. If $0\le u\le1$, the convolution remains in $[0,1]$ and has support in a fixed enlargement.
\end{lemma}

\begin{proof}
Right translation by a point in $B_\sigma$ changes $u$ in $L^1$ by at most $C\sigma|D_Hu|(G)$, as follows by integrating the directional variation along a horizontal path from the identity. Thus $w_\sigma=u*\rho_\sigma$ satisfies
\[
 \|w_\sigma-u\|_1\le C\sigma|D_Hu|(G),\qquad
 \|\nabla_{\rm Eucl}w_\sigma\|_1\le C\sigma^{-r}\|u\|_1.
\]
The measure-preserving left-translation argument in the proof of Lemma~\ref{lem:smoothing} now gives
\[
 \|\rho_\varepsilon*u-u\|_1
 \le2C\sigma|D_Hu|(G)+C\varepsilon\sigma^{-r}\|u\|_1.
\]
We take $\sigma=\varepsilon^{1/(r+1)}$. Left convolution is retained because it does not increase perimeter by \eqref{eq:convolution-perimeter}.

For the derivative estimate, we differentiate $\rho_\varepsilon(xy^{-1})$ in $x$. The coordinate coefficients are bounded on the fixed support. Changing variables by Haar invariance and using $\|D^j\rho_\varepsilon\|_1\le C\varepsilon^{-jr}$ then bounds integration against $u$ by $C\varepsilon^{-2r}\|u\|_\infty$. Positivity gives the remaining assertions.
\end{proof}

\begin{proposition}\label{prop:rate-upper}
With $\beta=\aexp\gamma/[Q(1+\aexp)]$, one has
$I(n)\le c_\phi n^{\aexp}+Cn^{\aexp-\beta}$.
\end{proposition}

\begin{proof}
Let $E$ be a bounded minimizer with $M\mu(E)=1$. Its dilates $E'$ with scale in $[1/2,1]$ lie in one compact set and have uniformly bounded volume and perimeter. We smooth $\one_{E'}$ by left convolution and apply Lemma~\ref{lem:quant-recovery}. The derivative estimate in Lemma~\ref{lem:continuum-smooth-error} gives
\[
 \E_h(f)\le\Per_\phi(E')+Ch\varepsilon^{-2r},\qquad
 \sum_\alpha\|f_h^\alpha-\one_{E'}\|_1
 \le C(\varepsilon^\theta+h\varepsilon^{-2r}).
\]
We take $\varepsilon=h^{1/(2r+\theta)}$, so both errors are $O(h^\gamma)$.

We balance the loss from restricting the levels against the boundary cost of adding missing vertices. We put $\tau=h^{\gamma/(1+\aexp)}$, so that $h^\gamma=\tau^{1+\aexp}$. By \eqref{eq:level-distance}, every level in $[\tau^{\aexp},1-\tau^{\aexp}]$ has normalized mass within $C_1\tau$ of $M\mu(E')$. We choose the dilate with $M\mu(E')=1-2C_1\tau$. The constants above are uniform in this choice. Coarea supplies a common level defining a finite set $A$ for which
\[
 h^{Q-1}\partial(A)\le\frac{\E_h(f)}{1-2\tau^{\aexp}}
 \le\Per_\phi(E')+C\tau^{\aexp}\le c_\phi+C\tau^{\aexp}.
\]
Also $1-3C_1\tau\le h^Q|A|\le1-C_1\tau$. For small $h$, the required dilating scale lies in $[1/2,1]$.

We set $h=n^{-1/Q}$. Then $s=n-|A|\le Cn\tau$. The coarse upper bound in Proposition~\ref{prop:sharp-law} gives $I(s)\le Cs^{\aexp}$. Adding a set of this size far from $A$ yields
\[
 I(n)\le c_\phi n^{\aexp}+Cn^{\aexp}\tau^{\aexp}
 =c_\phi n^{\aexp}+Cn^{\aexp-\beta}.
\]
We increase $C$ to cover the finitely many remaining cardinalities.
\end{proof}

Both estimates hold with the sharper exponent
\begin{equation}\label{eq:explicit-beta}
 \beta=\frac{\aexp}{Q(1+\aexp)(2r(r+1)+1)}
 =\frac{Q-1}{Q(2Q-1)(2r(r+1)+1)}.
\end{equation}
The first layer has dimension at least two, and every layer up to $r$ is nonzero. Hence $Q\ge1+r(r+1)/2$, so $2r(r+1)+1\le4Q-3$. Since $Q\ge2$,
\[
 \beta\ge\frac{Q-1}{Q(2Q-1)(4Q-3)}\ge\frac1{8Q^2}.
\]
Proposition~\ref{prop:concentration} completes Theorem~\ref{thm:main}. The power estimate requires no classification, boundary regularity, or quantitative stability of continuum isoperimetric regions.

The profile estimate also determines the approximate size of the smallest finite set whose boundary is at most $\varepsilon$ times its size. To obtain this F\o lner estimate, we invert the leading relation $I(n)/n\sim c_\phi n^{-1/Q}$ and control the error.

\begin{corollary}\label{cor:monotone-folner}
For either boundary in Theorem~\ref{thm:main}, we define
\[
 I_{\mathrm{up}}(n)=\inf_{m\ge n}I(m),\qquad
 F_\partial(\varepsilon)=\min\{|A|:0<|A|<\infty,\ \partial(A)\le\varepsilon|A|\}.
\]
With $\beta$ as in \eqref{eq:profile-rate}, one has
\[
 I_{\mathrm{up}}(n)=c_\phi n^{\aexp}+O(n^{\aexp-\beta}),\qquad
 \frac{I(n)}{I_{\mathrm{up}}(n)}=1+O(n^{-\beta}),
\]
and
\[
 F_\partial(\varepsilon)=\left(\frac{c_\phi}{\varepsilon}\right)^Q
 \bigl(1+O(\varepsilon^{Q\beta})\bigr)\qquad(\varepsilon\downarrow0).
\]
\end{corollary}

\begin{proof}
For sufficiently large $n$ and every $m\ge n$, Theorem~\ref{thm:main} gives $I(m)\ge(c_\phi-Cn^{-\beta})m^{\aexp}\ge c_\phi n^{\aexp}-Cn^{\aexp-\beta}$. Since $I_{\mathrm{up}}(n)\le I(n)$, the first two assertions follow.

We put $N=(c_\phi/\varepsilon)^Q$. The coarse bounds $I(k)\asymp k^{\aexp}$ first give $F_\partial(\varepsilon)\asymp N$. Uniformly for $k\asymp N$,
\[
 \frac{I(k)}k=\varepsilon\left(\frac Nk\right)^{1/Q}
 \bigl(1+O(N^{-\beta})\bigr).
\]
For a sufficiently large fixed $K$, this excludes every admissible $k$ in that range below $N(1-KN^{-\beta})$, whereas $\lceil N(1+KN^{-\beta})\rceil$ is admissible. Since $\beta<1$, integer rounding does not affect the error. These bounds prove the last assertion.
\end{proof}

\section{Virtually abelian groups}\label{sec:example}

When $G=\R^d$, the Wulff body $W$ minimizes $\Per_\phi$ among sets of volume $\mu(W)$. Its dilates give minimizers at every positive volume. The anisotropic Euclidean isoperimetric inequality states that
\begin{equation}\label{eq:wulff-constant}
 C_\phi=d\,\mu(W)^{1/d}\qquad(G=\R^d),
\end{equation}
with equality for dilates and translates of $W$; see \cite[Sections~9.2--9.3, Proposition~30]{PansuPeriodic} and \cite[Section~4.1]{DelNinPetrache}. The discrete leading constant is $c_\phi=M^{-(d-1)/d}C_\phi$, as in \eqref{eq:constant}.

For Theorem~\ref{thm:abelian-stability}, we apply the quantitative Wulff inequality to a level set of the smoothed indicators, then transfer the estimate to the discrete set.

\begin{proof}
We use lattice-basis coordinates, so $\mu=dx$. We put $h=n^{-1/d}$, $E_h=h^{d-1}\partial(A)$, $u^\alpha=u_h^\alpha$, and $\bar u=M^{-1}\sum_\alpha u^\alpha$. The path estimate gives, for every $y\in\R^d$,
\begin{equation}\label{eq:abelian-cell-translations}
 \|u^\alpha(\,\cdot+y)-u^\alpha\|_1\le C(|y|+h)E_h,\qquad
 \|u^\alpha-\bar u\|_1\le ChE_h.
\end{equation}
Indeed, we put $c_\kappa=\mu(D\cap(\kappa+D-y/h))$. These nonnegative coefficients sum to one, and $c_\kappa>0$ implies $|\kappa|_S\le C(1+|y|/h)$ because $D$ is bounded. Intersecting the translated cells expresses the first norm as $h^d\sum_\kappa c_\kappa\sum_\lambda |f_{\lambda+\kappa,\alpha}-f_{\lambda,\alpha}|$, so Lemma~\ref{lem:paths} applies. The second assertion follows from \eqref{eq:coset-sync}. No perimeter assumption on $D$ is used.

We choose an even nonnegative smooth mollifier $\rho_\varepsilon$, and put $w=\rho_\varepsilon*\bar u$. Then $0\le w\le1$, $\int w=1/M$, and
\begin{equation}\label{eq:abelian-global-comparison}
 \|w-u^\alpha\|_1\le C(\varepsilon+h)E_h,\qquad
 \Per_\phi(w)\le(1+Ch/\varepsilon)E_h.
\end{equation}
The first estimate follows from \eqref{eq:abelian-cell-translations}. To prove the second with its sharp leading coefficient, we fix $F\in C_c^\infty(\R^d;\Fcal)$ and let $F_\varepsilon=\rho_\varepsilon*F$. Its values remain in $\Fcal$, and $\|\nabla F_\varepsilon\|_\infty\le C/\varepsilon$. Since $v_e\in\Lambda$ and translation by $hv_e$ carries each cell to the corresponding cell, pointwise duality for the discrete energy gives
\[
 J=\frac1h\sum_e\int F_{\varepsilon,e}(x)
       \bigl(u^{\beta(e)}(x+hv_e)-u^{\alpha(e)}(x)\bigr)\,dx\le E_h.
\]
After changing variables in the first term, the zeroth-order terms cancel by $\Div_0F_\varepsilon=0$. We write $\partial_{v_e}$ for the Euclidean directional derivative along $v_e$. For each remaining term, the error between its difference quotient and $-\partial_{v_e}F_{\varepsilon,e}$, integrated against $u^{\beta(e)}$, is
\[
 \int_0^1\int
 \bigl(u^{\beta(e)}(x)-u^{\beta(e)}(x+thv_e)\bigr)
 \partial_{v_e}F_{\varepsilon,e}(x)\,dx\,dt.
\]
Its absolute value is at most $Ch\varepsilon^{-1}E_h$ by \eqref{eq:abelian-cell-translations}. Replacing $u^{\beta(e)}$ by $\bar u$ in the first-order terms has the same bound, by \eqref{eq:coset-sync}. Consequently
\[
 -\int w\sum_e\partial_{v_e}F_e
 =-\int\bar u\sum_e\partial_{v_e}F_{\varepsilon,e}
 \le E_h+Ch\varepsilon^{-1}E_h.
\]
Taking the supremum proves \eqref{eq:abelian-global-comparison} for both boundaries.

It suffices to consider $E_h\le2c_\phi$ and sufficiently small $h$; the other cases follow by enlarging $C$. We set $\varepsilon=h^{1/2}$ and $\eta=Ch^{1/2}$. If $v_0=\int u^1$, then $v_0=1/M+O(h)$ by \eqref{eq:abelian-cell-translations}. For $E_t=\{w>t\}$ and $v_t=\mu(E_t)$, the layer-cake formula gives
\[
 \int_0^1|v_t-v_0|\,dt\le\|w-u^1\|_1\le\eta.
\]
We put $\aexp=(d-1)/d$. The elementary inequality $v^\aexp\ge v_0^\aexp-v_0^{\aexp-1}(v_0-v)_+$ implies
\[
 \int_0^1v_t^\aexp\,dt\ge M^{-\aexp}-C\eta.
\]
Since the isoperimetric deficit $\Delta_t=\Per_\phi(E_t)-C_\phi v_t^\aexp$ is nonnegative, coarea and \eqref{eq:abelian-global-comparison} yield
\[
 \int_0^1\Delta_t\,dt\le E_h-c_\phi+C\eta
 \le C\bigl(\delta_\phi(A)+h^{1/2}\bigr).
\]
We choose $t\in[1/4,3/4]$ with $\Delta_t\le C(\delta_\phi(A)+h^{1/2})$. For every such level, $\|\one_{E_t}-u^1\|_1\le4\eta$; hence $v_t=1/M+O(\eta)$ and $v_t$ is bounded away from zero. The quantitative Wulff inequality gives a translate of $W_{v_t}=(v_t/\mu(W))^{1/d}W$ such that
\[
 \mu\bigl(E_t\triangle(x+W_{v_t})\bigr)
 \le C\bigl(\delta_\phi(A)+h^{1/2}\bigr)^{1/2}.
\]
The concentric dilates $W_{v_t}$ and $W_*$ differ in measure by $|v_t-1/M|$. Combining this with the level estimate and \eqref{eq:coset-sync} gives
\begin{equation}\label{eq:abelian-continuum-stability}
 \sum_\alpha\|u^\alpha-\one_{x+W_*}\|_1
 \le C\bigl(\delta_\phi(A)+h^{1/2}\bigr)^{1/2}.
\end{equation}
Finally we choose $\kappa\in\Lambda$ with $|h\kappa-x|\le Ch$. Translating the fixed convex body by this amount changes its measure by $O(h)$. Its approximation by cells with associated lattice points in $h\kappa+W_*$ also differs from it by $O(h)$: every differing cell lies within distance $Ch$ of its boundary, whose $Ch$-neighborhood has measure $O(h)$. This approximation is exactly the cell function of $\kappa\mathcal W_n$ in every coset. Thus \eqref{eq:abelian-continuum-stability}, multiplied by $h^{-d}=n$, proves \eqref{eq:abelian-stability}.
\end{proof}

We next compare the large-scale word metric with the limiting vertex perimeter. Let $\Lambda\lhd\Gamma$ be free abelian of rank $d\ge2$ and index $M$. We put $V=\Lambda\otimes_{\Z}\R$. Conjugation induces a finite group $H\le\mathrm{GL}(V)$. In additive notation on $\Lambda$, the limit $\|\lambda\|_S=\lim_{k\to\infty}k^{-1}|k\lambda|_S$ extends to a norm on $V$, called the stable norm. We write $B_S$ for its unit ball and $W_{v,S}$ for the Wulff body with support function $\phi_{v,S}$. The measured asymptotic cone is $V$ with this norm and measure $M\mu$. If perimeter is defined relative to that measure rather than $\mu$, its support body is $W_{v,S}/M$.

\begin{theorem}\label{thm:holonomy-characterization}
The following are equivalent:
\begin{enumerate}[(i)]
 \item $H\subset\{I,-I\}$;
 \item $W_{v,S}=MB_S$ for every finite symmetric generating set $S$;
 \item $c_{v,S}=c_{v,T}$ whenever $S,T$ are finite symmetric generating sets whose word lengths differ by a bounded amount.
\end{enumerate}
If (i) fails, there are $S\subset T$ with $d_T\le d_S\le d_T+C$ and $c_{v,S}<c_{v,T}$.
\end{theorem}

\begin{proof}
We first rewrite the cell formula using the finite conjugation action $H$ instead of all cosets in $\Gamma/\Lambda$. We choose $q_1=e$ and let $\rho:\Gamma/\Lambda\to H$ be the conjugation map. For $\gamma\in\Gamma$, we also write $\rho(\gamma)$ for the image of its class in $\Gamma/\Lambda$. In additive notation on $\Lambda$, we set
\[
 c(f,g)=q_fq_gq_{fg}^{-1},\qquad x_f=M^{-1}\sum_{k\in\Gamma/\Lambda}c(f,k).
\]
Here $c(f,g)\in\Lambda$ records the failure of the representatives to multiply exactly. Associativity gives $c(f,g)+c(fg,k)=\rho(f)c(g,k)+c(f,gk)$. Averaging over $k$ yields $c(f,g)=x_f+\rho(f)x_g-x_{fg}$. Thus the map $b(\lambda q_f)=\lambda+x_f\in V$ satisfies $b(\gamma\gamma')=b(\gamma)+\rho(\gamma)b(\gamma')$. If $\bar s\in\Gamma/\Lambda$ represents $s$, then $v_{fs}=x_f+\rho(f)b(s)-x_{f\bar s}$. We substitute $\chi_f=\psi_f-p(x_f)$ in \eqref{eq:vertex-cell}. The resulting convex objective is invariant under $\chi_f\mapsto\chi_{nf}$ for $n\in\ker\rho$. Averaging makes $\chi$ constant on fibers of $\Gamma/\Lambda\to H$ without increasing the objective. Thus, with $m=|H|$,
\begin{equation}\label{eq:holonomy-reduced-cell}
 \phi_{v,S}(p)=\frac{M}{m}\min_{\eta\in\R^H}\sum_{\omega\in H}
 \max\left(0,\max_{s\in S}\{p(\omega b(s))+\eta_{\omega\rho(s)}-\eta_\omega\}\right).
\end{equation}
The coordinates $b(\lambda q_f)$ differ from $\lambda$ by a bounded amount.

If $H=\{I\}$, formula~\eqref{eq:holonomy-reduced-cell} gives $W_{v,S}/M=\conv b(S)$, and this is $B_S$. For $H=\{I,-I\}$, we put
\[
 P=\{b(s):s\in S,\ \rho(s)=I\},\qquad J=\{b(s):s\in S,\ \rho(s)=-I\}.
\]
Here $J\ne\varnothing$ and $P=-P$. For fixed $p$, let $a=\max(0,\max_{P}p)$, taking $a=0$ if $P$ is empty, and $a_\pm=\max_{J}p,\min_{J}p$, respectively. The reduced minimum equals
\[
 \min_t\{\max(a,a_++t)+\max(a,-a_--t)\}
 =\max(2a,a_+-a_-),
\]
with equality at $t=-(a_++a_-)/2$. Hence $W_{v,S}/M=\conv(P\cup(J-J)/2)$, inserting $0$ if $P$ is empty. This is $B_S$: successive letters with $\rho(s)=-I$ contribute pairs $u-u'$, $u,u'\in J$, while the other letters contribute elements of $P$; at most one bounded term remains. Conversely, powers of generators with $\rho(s)=I$ and of products of two generators with $\rho(s)=-I$ realize these velocities, and concatenation realizes their convex hull. This proves (i)$\Rightarrow$(ii). Bounded additive distance gives equal stable balls, so \eqref{eq:wulff-constant} proves (ii)$\Rightarrow$(iii).

We now suppose $H\not\subset\{I,-I\}$. We construct two generating sets with the same word metric up to a bounded additive error but different leading vertex constants. We choose a finite symmetric $R\subset\Gamma\setminus\Lambda$ whose images generate $\Gamma/\Lambda$. There is a symmetric lattice polytope $K$ containing $b(R)$ and a lattice basis but not invariant under $H$: we choose $\omega_0\in H\setminus\{I,-I\}$ and $v\in\Lambda$ with $\omega_0v\notin\R v$. We start with a symmetric lattice polytope containing these prescribed points and enlarge it by $\pm Lv$. Its projection to $V/\R v$ stays bounded as $L\to\infty$, while the projection of $L\omega_0v$ does not. We choose a finite symmetric $P\subset\Lambda\setminus\{0\}$ generating $\Lambda$ with $\conv P=K$, and put
\[
 S=P\cup R,\qquad T=HP\cup R,\qquad K_H=\conv\bigcup_{\omega\in H}\omega K.
\]
At coset $\omega$, translation generators begin and end in the same coset, so their shift terms cancel. Their contribution is at least $h_K(\omega^*p)$, where $h_K(p)=\max_{x\in K}p(x)$ is the support function of $K$ and $(\omega^*p)(x)=p(\omega x)$. Since $b(R)\subset K$, taking $\eta=0$ in \eqref{eq:holonomy-reduced-cell} attains their sum. Therefore
\begin{equation}\label{eq:holonomy-orbit-formula}
 W_{v,S}=\frac{M}{m}\sum_{\omega\in H}\omega K,\qquad W_{v,T}=MK_H.
\end{equation}
Every affine increment belongs to $K_H$, and conjugated translation powers realize $HP$, so $B_S=B_T=K_H$. Each $\omega K\subset K_H$. Equality of their Minkowski average with $K_H$ would force every support function $h_K(\omega^*p)\le h_{K_H}(p)$ to equal $h_{K_H}(p)$, implying $K=K_H$. Thus $W_{v,S}\subsetneq W_{v,T}$, and strict volume inequality gives $c_{v,S}<c_{v,T}$.

It remains to bound the difference of the word metrics. We fix $S$-words $u_\omega$ with $\rho(u_\omega)=\omega$. We move the translation letters of a $T$-word to the left, retaining its $R$-subword in its original order, and group the resulting translations $\omega a$ according to $\omega$. Each group is represented by $u_{\omega}w_{\omega}u_{\omega}^{-1}$ with the same number of $P$-letters. These conjugates lie in the abelian subgroup $\Lambda$; their product followed by the retained $R$-subword represents the original element. Hence $d_T\le d_S\le d_T+2\sum_\omega|u_\omega|_S$, contradicting (iii).
\end{proof}

\begin{proposition}\label{prop:holonomy-orbit}
Let $\Gamma=\Lambda\rtimes H$, let $P=-P\subset\Lambda\setminus\{0\}$ be finite with $HP$ generating $\Lambda$, and let $R=R^{-1}\subset H\setminus\{1\}$ generate $H$. We set
\[
 S=(P\times\{1\})\cup(\{0\}\times R),\qquad
 T=(HP\times\{1\})\cup(\{0\}\times R),\qquad K=\conv P.
\]
Then $B_S=B_T=\conv\bigcup_\omega \omega K$, $W_{v,S}=\sum_\omega \omega K$, $W_{v,T}=|H|B_S$, and $d_T\le d_S\le d_T+2(|H|-1)$. The two constants agree exactly when $K$ is $H$-invariant.
\end{proposition}

\begin{proof}
The translation-generator calculation in \eqref{eq:holonomy-orbit-formula} applies here. The $R$-edges have zero displacement, so they add no term when all shifts are zero. For the metric bound, we group translation letters as above and insert the translation letters assigned to $\omega$ while visiting $\omega$ on a closed walk through the quotient. Traversing a spanning tree twice uses $2(|H|-1)$ letters. We append the original $R$-subword. The support-function argument above proves the equality criterion.
\end{proof}

We write $\mathcal L^d$ for Lebesgue measure on $\R^d$. For example, we take $H=\langle\tau\rangle$ exchanging the coordinates, $P=\{\pm e_1\}$, and $R=\{\tau\}$. Then $d_T\le d_S\le d_T+2$, the common measured cone is the $\ell^1$ plane with measure $2\mathcal L^2$, and the Wulff bodies have areas $4$ and $8$. Thus $c_{v,S}=2\sqrt2$ and $c_{v,T}=4$. The two limits have the same volume--diameter problem, since their metric and measure agree, but different perimeter problems: the vertex Wulff bodies are a square and a diamond, respectively. This separates the isodiametric question from the isoperimetric one.

\begin{corollary}\label{cor:holonomy-unbounded}
For every $d\ge2$ there are finite symmetric generating sets $S_d\subset T_d$ of a virtually abelian group of rank $d$, with bounded additive metric difference and common measured cone $(\R^d,\|\cdot\|_\infty,2\mathcal L^d)$. Let $\mathsf Z_d$ be a $d\times d$ random matrix whose entries are independent and equal to $1$ or $-1$ with probability $1/2$, and let $\mathbb E$ denote expectation. Their leading constants satisfy
\begin{equation}\label{eq:holonomy-gap}
 \frac{c_{v,T_d}}{c_{v,S_d}}
 =\left(\frac{d!}{\mathbb E|\det\mathsf Z_d|}\right)^{1/d}
 \ge(d!)^{1/(2d)}\sim\sqrt{d/e},
\end{equation}
where the asymptotic equivalence is as $d\to\infty$.
\end{corollary}

\begin{proof}
Let $H_d$ be the coordinate sign-change group, $v=(1,\ldots,1)$, and
\[
 \Lambda_d=\operatorname{span}_{\Z}\{\varepsilon:\varepsilon\in\{-1,1\}^d\}
 =\{x\in\Z^d:x_i\equiv x_j\pmod2\text{ for all }i,j\}.
\]
We apply Proposition~\ref{prop:holonomy-orbit} with $P=\{v,-v\}$ and $R$ the coordinate sign changes. Then $M=2^d$, $K_H=[-1,1]^d$, and
\[
 W_{v,S_d}=\sum_{\varepsilon\in\{-1,1\}^d}[-\varepsilon,\varepsilon],
 \qquad W_{v,T_d}=M[-1,1]^d.
\]
The lattice has covolume $2^{d-1}$, so $M\mu=2\mathcal L^d$. The zonotope volume formula, averaged over ordered $d$-tuples of sign vectors, gives
\[
 |W_{v,S_d}|=\frac{2^dM^d}{d!}\mathbb E|\det\mathsf Z_d|,
 \qquad |W_{v,T_d}|=2^dM^d.
\]
The common normalization cancels in the ratio of constants. Any two distinct permutation terms in the determinant expansion have expected product zero, so $\mathbb E(\det\mathsf Z_d)^2=d!$. Cauchy--Schwarz and Stirling's formula finish the proof.
\end{proof}

For Riemannian metrics lifted from tori, Burago--Ivanov \cite[p.~787]{BuragoIvanov} show that the stable norm and asymptotic volume determine the isoperimetric constant in dimension two, but not in higher dimensions. The finite-quotient example above shows that these data need not determine the isoperimetric constant even in rank two. Breuillard and Le Donne \cite{BreuillardLeDonne} give quantitative convergence estimates for rescaled nilpotent word-metric balls. In the two-dimensional example above, the projection to $\Z^2$ has bounded additive error for both metrics, hence error $O(R^{-1})$ after rescaling by $R^{-1}$. Even this rate of convergence to a common measured cone does not determine the limiting perimeter.

\begin{proposition}\label{prop:wallpaper-constants}
For the generating sets below, the edge weights are one. The wallpaper groups $pg,p3,p4$ then have the following constants. We use a translation-lattice basis of covolume one and dual covector coordinates $p_1,p_2$. We put $\phi_\triangle(p)=|p_1|+|p_2|+|p_1+p_2|$.
\[
\begin{array}{c|c|c|c|c|c}
 \Gamma&M&\phi_v(p)&c_v&\phi_e(p)&c_e\\ \hline
 pg&2&\max(|p_1|,2|p_2|)&2\sqrt2&|p_1|+2|p_2|&4\\
 p3&3&\phi_\triangle(p)&4&\phi_\triangle(p)&4\\
 p4&4&|p_1|+|p_2|&2&|p_1|+|p_2|&2
\end{array}
\]
\end{proposition}

\begin{proof}
For $pg$, we take $g(u,v)=(u+1/2,-v)$, $b(u,v)=(u,v+1)$, and $S=\{g^{\pm1},b^{\pm1}\}$. The translation subgroup is $\langle g^2,b\rangle$, with representatives $e,g$. Writing $d=\psi_2-\psi_1$, the cell formulas give
\begin{align*}
 \phi_v(p)&=\min_d\bigl[\max(|p_2|,d,d-p_1)+\max(|p_2|,p_1-d,-d)\bigr]=\max(|p_1|,2|p_2|),\\
 \phi_e(p)&=2|p_2|+\min_d(|d|+|p_1-d|)=|p_1|+2|p_2|.
\end{align*}
The first sum is bounded below by $2|p_2|$ and $\pm p_1$, with equality at $d=p_1/2$. The vertex Wulff body has vertices $(\pm1,0),(0,\pm2)$ and area four; the edge Wulff body is a rectangle of area eight. At $p=(4,1)$, taking $\psi=0$ gives vertex energy five instead of four.

For $p3$, $a$ is the translation by $e_1$ in a triangular basis and $\tau$ satisfies $\tau e_1=e_2$, $\tau e_2=-e_1-e_2$. We take $S=\{a^{\pm1},\tau^{\pm1}\}$. Translation loops in the three cosets have directions $e_1,e_2,-e_1-e_2$, giving the lower bound $\phi_\triangle(p)$ in both cell formulas; $\psi=0$ attains it. The corresponding zonotope has area twelve.

For $p4$, we take $\tau(u,v)=(-v,u)$, $h(u,v)=(1-u,-v)$, and $S=\{\tau,\tau^{-1},h\}$. With representatives $\tau^i$, the rotation edges have zero displacement and the $h$ edge from $i$ to $i+2$ has displacement $\tau^ie_1$. Taking $\psi=0$ bounds both integrands above by $|p_1|+|p_2|$. Each choice of two $h$ edges with displacements $\pm e_1,\pm e_2$ extends by two rotation edges to a directed four-cycle through all cosets. The edge counts of these cycles are admissible in both dual problems, and their images under $T$ are the four vertices of $[-1,1]^2$. Thus both Wulff bodies have area four. Finally, \eqref{eq:constant} and \eqref{eq:wulff-constant} give $c_\phi=2\sqrt{|W|/M}$ throughout.
\end{proof}

The restriction $Q>1$ is necessary for an exact profile limit. On the unweighted ladder $\Z\times C_2$, where $C_2$ is the cyclic group of order two, we take the two horizontal generators and the vertical involution. Then the edge profile satisfies $I_e(2n)=4$ and $I_e(2n+1)=5$ for $n\ge1$. If both rows are nonempty their horizontal boundaries contribute at least four; odd cardinality forces a vertical crossing. If only one row is nonempty the boundary is at least its cardinality plus two. Full columns, with one additional vertex when needed, give equality. This dependence on parity is no longer negligible when the isoperimetric exponent is zero.

\section{The Heisenberg group}\label{sec:heisenberg}

We now bound the limiting perimeters in the standard Heisenberg lattice and construct finite sets of every cardinality with an explicit upper bound. We use exponential coordinates $(x,y,z)$ with group law
\[
 (x,y,z)(x',y',z')=(x+x',y+y',z+z'+\tfrac12(xy'-yx')),
\]
with horizontal vector fields $X=\partial_x-\frac y2\partial_z$ and $Y=\partial_y+\frac x2\partial_z$, and Haar measure $dx\,dy\,dz$. For a region $E$ in the continuous group, we write $|E|$ for this volume. The lattice $x,y,z+xy/2\in\Z$ has covolume one. The standard generators are $(\pm1,0,0)$ and $(0,\pm1,0)$, with unit edge weights. For a horizontal covector $p$, let $p_1,p_2$ be its components in the basis dual to $X,Y$. The inverse pairs give $\phi_e(p)=|p_1|+|p_2|$ and $\phi_v(p)=\max(|p_1|,|p_2|)$; we write $c_e,c_v$ for their constants.

\begin{proposition}\label{prop:heisenberg-bounds}
For the standard Heisenberg generating set,
\[
 \frac4{3^{1/4}}\le c_e\le\frac8{3^{3/4}},\qquad c_v=\frac{c_e}{\sqrt2}.
\]
Theorem~\ref{thm:main} holds with $\beta=3/364$.
\end{proposition}

\begin{proof}
For the lower bound, we estimate volume by the areas of two horizontal projections, then bound those areas by perimeter; see \cite[Sections~2 and~4]{FasslerPinamonti}. For $B\subset\R^2$ of finite area, we put $T\one_B(s,t)=\int\one_B(x,t+sx)\,dx$. Expanding the cube, ordering $x_1<x_2<x_3$, and changing $(s,t)$ to $z_i=t+sx_i$, $i=1,3$, gives
\[
 \int(T\one_B)^3=6\int_{x_1<x_3}\one_B(x_1,z_1)\one_B(x_3,z_3)J\,d(x_1,z_1,x_3,z_3)\le3|B|^2,
\]
where $J=(x_3-x_1)^{-1}\int_{x_1}^{x_3}\one_B(x,z_1+(x-x_1)(z_3-z_1)/(x_3-x_1))\,dx\le1$. Let $E$ have finite volume and perimeter. We define its essential horizontal projections $\Pi_X,\Pi_Y$ to be the parameters of the $X$- and $Y$-fibers whose intersections with $E$ have positive one-dimensional measure. In coordinates $(y,z+xy/2)$ and $(x,z-xy/2)$, we write $a,b$ for their Lebesgue areas. Slicing gives $|D_X\one_E|\ge2a$ and $|D_Y\one_E|\ge2b$. With $k=z+xy/2$, H\"older's inequality yields
\[
 |E|\le\int\one_{\Pi_X}(y,k)\one_{\Pi_Y}(x,k-xy)\,dx\,dy\,dk
 \le a^{2/3}\|T\one_{\Pi_Y}\|_3\le3^{1/3}(ab)^{2/3}.
\]
Hence $\Per_{\phi_e}(E)\ge2(a+b)\ge4\sqrt{ab}\ge4\,3^{-1/4}|E|^{3/4}$.

For the upper bound, we consider the square-based region
\begin{equation}\label{eq:square-bubble}
 E_\square=\{(x,y,z): |x|,|y|\le1,\ 2|z|\le1-|xy|\}.
\end{equation}
It has volume three. After central translation by $1/2$, its bounds in $k=z+xy/2$ are $\max(xy,0)\le k\le1+\min(xy,0)$. Both horizontal projections are rectangles of area two, and every nonempty horizontal fiber is an interval. Thus $\Per_{\phi_e}(E_\square)=8$.

After dilation and central translation, this is the square case of the construction in \cite[Section~8.1, pp.~34--35]{FranceschiMontiRighiniSigalotti}. Their horizontal vector norm in this case is $\ell^\infty$, whose dual is the perimeter density $\phi_e=\ell^1$ used here. The limiting square norms are also discussed in \cite[Remark~4.7]{PozueloRitore}. We use the region only as a competitor, without assuming that it minimizes perimeter among all sets.

The planar rotation $R_{\pi/4}$ through $\pi/4$ extends to a measure-preserving graded automorphism $\mathcal R$. Since $\|R_{\pi/4}p\|_1=\sqrt2\|p\|_\infty$, one has $\Per_{\phi_e}(\mathcal R E)=\sqrt2\Per_{\phi_v}(E)$. Taking infima proves the relation between the constants. Equation~\eqref{eq:explicit-beta}, with $Q=4,r=2$, gives the exponent.
\end{proof}

In integer matrix coordinates the group law is $(i,j,k)(i',j',k')=(i+i',j+j',k+k'+ij')$. The standard generator edges join $(i,j,k)$ to $(i\pm1,j,k)$ and $(i,j\pm1,k\pm i)$.

\begin{proposition}\label{prop:heisenberg-exact-upper}
For the standard Heisenberg generating set, there is an absolute constant $C$ such that
\[
 I_e(n)\le\frac8{3^{3/4}}n^{3/4}+Cn^{1/4}\qquad(n\ge1).
\]
\end{proposition}

\begin{proof}
We first take lattice approximations to dilates of $E_\square$, allowing the height to vary. We then add the remaining vertices along generator paths in the top layer. For integers $N\ge1$ and $K\ge N^2$, we put
\[
 A_{N,K}=\{(i,j,k):|i|,|j|\le N,\ \max(ij,0)\le k\le K+\min(ij,0)\}.
\]
Summing the vertical intervals and counting the interval fibers in the two generator directions gives
\[
 v_{N,K}:=|A_{N,K}|=(2N+1)^2(K+1)-N^2(N+1)^2,\qquad
 B(A_{N,K})=4(2N+1)(K+1).
\]
In particular, $A_N=A_{N,N^2}$ satisfies $|A_N|=3N^4+2N^3+4N^2+4N+1$ and $B(A_N)=4(2N+1)(N^2+1)$.

For large $n$, we choose the largest $N$ with $|A_N|\le n$, and set
\[
 K=N^2+\left\lfloor\frac{n-|A_N|}{(2N+1)^2}\right\rfloor.
\]
Then $n\asymp N^4$, $K-N^2=O(N)$, and $0\le n-v_{N,K}<(2N+1)^2$. Writing $K=N^2+t$ and $C_\square=8/3^{3/4}$, Taylor's formula gives
\[
 C_\square v_{N,K}^{3/4}=8N^3+4N^2+8Nt+O(N)=B(A_{N,K})+O(N).
\]

We add the remaining vertices from the top layer
\[
 L_{N,K}=A_{N,K+1}\setminus A_{N,K}
 =\{(i,j,K+1+\min(ij,0)):|i|,|j|\le N\}.
\]
We begin with $(0,0,K+1)$, then add the four arms of the cross $ij=0$ outward from the origin. The first point increases the edge boundary by at most four; each later point has its predecessor already present and increases it by at most two. Each initial segment increases the boundary by at most $4+2(4N)=8N+4$. The other points form $4N$ chains of length $N$: we use consecutive $i$ with $j$ fixed in the quadrants $ij>0$, and consecutive $j$ with $i$ fixed in the quadrants $ij<0$. These are generator chains. Each point $(i,j,k)$ on them has a neighbor in $A_{N,K}$, namely $(i,j-\operatorname{sgn}i,k-|i|)$ when $ij>0$, and $(i+\operatorname{sgn}j,j,k)$ when $ij<0$. The defining inequalities follow from $K\ge N^2$.

We append these chains one at a time, each oriented outward from its endpoint on the cross. Every new point then has two distinct neighbors already present: its predecessor on the chain or cross, and the neighbor in $A_{N,K}$ just described. Since the graph has degree four, adding such a point changes the boundary by $4-2m\le0$, where $m\ge2$ is the number of its neighbors already present. Thus every initial segment of the full ordering increases the boundary by at most $8N+4$. Taking the initial segment of cardinality $n-v_{N,K}$ gives a set of size $n$ with boundary at most $C_\square v_{N,K}^{3/4}+O(N)\le C_\square n^{3/4}+O(n^{1/4})$. Increasing $C$ covers the remaining values of $n$.
\end{proof}

Weights $a,b>0$ on the two inverse pairs give $c_e(a,b)=\sqrt{ab}\,c_e(1,1)$ by a volume-preserving diagonal graded automorphism.

For the standard word balls $\mathcal B_R=B_S(e,R)\subset\Gamma$, the outer vertex boundary is exactly $\mathcal B_{R+1}\setminus \mathcal B_R$. Shapiro's sphere formula, recorded in \cite[Section~10.4.2, p.~262]{DuchinShapiro}, gives $|\mathcal B_R|=(31/72)R^4+O(R^3)$ and $V(\mathcal B_R)=(31/18)R^3+O(R^2)$. Theorem~\ref{thm:main} and Proposition~\ref{prop:heisenberg-bounds} therefore imply
\begin{equation}\label{eq:heisenberg-ball-deficit}
 \lim_{R\to\infty}\frac{V(\mathcal B_R)}{I_v(|\mathcal B_R|)}
 =\frac{4(31/72)^{1/4}}{c_v}
 \in\left[\left(\frac{93}{32}\right)^{1/4},\left(\frac{31}{6}\right)^{1/4}\right].
\end{equation}
In particular, these balls have asymptotic vertex boundary at least $30.56\%$ larger than the minimum.

For regions between two Lipschitz graphs over the fixed square, we compare the graph gradients with explicit bounded vector fields. Integration by parts bounds the perimeter, and the equality conditions identify the minimizer.

\begin{proposition}\label{prop:heisenberg-square-graphs}
Let $f,g$ be Lipschitz on $[-1,1]^2$, with $f<g$ in $\Omega=(-1,1)^2$, and $E=\{(x,y,z):(x,y)\in\Omega,\ f(x,y)<z<g(x,y)\}$. Then
\[
 \Per_{\phi_e}(E)\ge2|E|+2\ge\frac8{3^{3/4}}|E|^{3/4}.
\]
The unique minimizer of the perimeter quotient in this class is $E_\square$, up to a null set and a central translation.
\end{proposition}

\begin{proof}
We put $\mathbf a=(y/2,-x/2)$. To estimate the two graph terms, we use
\[
 \Theta=(\operatorname{sgn}(y)(1-|x|),-\operatorname{sgn}(x)(1-|y|))
\]
and the vector fields $\sigma_g=\Theta-(x,y)$, $\sigma_f=\Theta+(x,y)$. Both fields lie in $[-1,1]^2$. The jumps of $\Theta$ across each axis are tangent to that axis, so they contribute no distributional divergence. Direct calculation away from the axes gives $\Div \Theta=0$. Hence $\Div\sigma_g=-2$ and $\Div\sigma_f=2$ distributionally. On each side of $\partial\Omega$, their outward normal components are $-1$ and $1$. With $\mathcal H^1$ denoting arclength on $\partial\Omega$, the two graphs and vertical sides contribute
\begin{align*}
 \Per_{\phi_e}(E)
 &=\int_\Omega\bigl(\|\nabla g+\mathbf a\|_1+\|\nabla f+\mathbf a\|_1\bigr)
   +\int_{\partial\Omega}(g-f)\,d\mathcal H^1\\
 &\ge2\int_\Omega(g-f)+2\int_\Omega \Theta\cdot \mathbf a=2|E|+2.
\end{align*}
The integration by parts is performed in the four quadrants. The terms on the axes cancel because the jumps of $\Theta$ are tangent to them; the outer boundary terms cancel the contribution from the vertical sides. Also $2\Theta\cdot \mathbf a=|x|+|y|-2|xy|$, whose integral is two. For $V=|E|$, the quotient $(2V+2)V^{-3/4}$ is minimized at $V=3$.

For equality, we first consider $0<x,y<1$. Here
$\sigma_g=(1-2x,-1)$ and $\sigma_f=(1,-1+2y)$. Equality in the dual $\ell^1$ bounds forces $g_x+y/2=0$ and $f_y-x/2=0$ almost everywhere, since the corresponding field components have absolute value strictly less than one. Thus $g=-xy/2+a(y)$ and $f=xy/2+b(x)$ in that quadrant. The same calculation gives $g_x=-y/2$, $f_y=x/2$ in the third quadrant, and $g_y=x/2$, $f_x=-y/2$ in the second and fourth. The equations in the adjacent quadrants make $g$ constant along $x=0$ and $f$ constant along $y=0$. Continuity therefore makes $a$ and $b$ constant; matching the four quadrants gives $g=c_+-|xy|/2$ and $f=c_-+|xy|/2$. Volume three gives $c_+-c_-=1$, hence the claimed set up to a null set and a central translation. For related restricted minimization by calibration, see \cite[Theorem~3.1]{RitoreCalibration} and \cite[Theorem~6.3]{PozueloRitore}.
\end{proof}

\subsection*{Acknowledgements} The author thanks IIT Bombay for providing ideal working conditions.


\begin{thebibliography}{99}

\bibitem{AmbrosioGhezziMagnani}
L. Ambrosio, R. Ghezzi, and V. Magnani, \emph{BV functions and sets of finite perimeter in sub-Riemannian manifolds}, Ann. Inst. H. Poincar\'e C Anal. Non Lin\'eaire \textbf{32} (2015), 489--517.

\bibitem{BarberErde}
B. Barber and J. Erde, \emph{Isoperimetry in integer lattices}, Discrete Anal. (2018), Paper No.~7, 16 pp.

\bibitem{BarberErdeKeevashRoberts}
B. Barber, J. Erde, P. Keevash, and A. Roberts, \emph{Isoperimetric stability in lattices}, Proc. Amer. Math. Soc. \textbf{151} (2023), 5021--5029.

\bibitem{Breuillard}
E. Breuillard, \emph{Geometry of locally compact groups of polynomial growth and shape of large balls}, Groups Geom. Dyn. \textbf{8} (2014), 669--732.

\bibitem{BreuillardLeDonne}
E. Breuillard and E. Le Donne, \emph{On the rate of convergence to the asymptotic cone for nilpotent groups and subFinsler geometry}, Proc. Natl. Acad. Sci. USA \textbf{110} (2013), 19220--19226.

\bibitem{BuragoIvanov}
D. Burago and S. Ivanov, \emph{On asymptotic isoperimetric constant of tori}, Geom. Funct. Anal. \textbf{8} (1998), 783--787.

\bibitem{CarbottiDonPallaraPinamonti}
A. Carbotti, S. Don, D. Pallara, and A. Pinamonti, \emph{Local minimizers and gamma-convergence for nonlocal perimeters in Carnot groups}, ESAIM Control Optim. Calc. Var. \textbf{27} (2021), suppl., S11.

\bibitem{ChambolleGiacominiLussardi}
A. Chambolle, A. Giacomini, and L. Lussardi, \emph{Continuous limits of discrete perimeters}, ESAIM Math. Model. Numer. Anal. \textbf{44} (2010), 207--230.

\bibitem{ChambolleKreutz}
A. Chambolle and L. Kreutz, \emph{Crystallinity of the homogenized energy density of periodic lattice systems}, preprint, \href{https://arxiv.org/abs/2106.08111v4}{arXiv:2106.08111v4}.

\bibitem{CicaleseLeonardi}
M. Cicalese and G. P. Leonardi, \emph{Maximal fluctuations on periodic lattices: an approach via quantitative Wulff inequalities}, Comm. Math. Phys. \textbf{375} (2020), 1931--1944.

\bibitem{DelNinPetrache}
G. Del Nin and M. Petrache, \emph{Continuum limits of discrete isoperimetric problems and Wulff shapes in lattices and quasicrystal tilings}, Calc. Var. Partial Differential Equations \textbf{61} (2022), Article 226.

\bibitem{DuchinShapiro}
M. Duchin and M. Shapiro, \emph{The Heisenberg group is pan-rational}, Adv. Math. \textbf{346} (2019), 219--263.

\bibitem{FasslerPinamonti}
K. F\"assler and A. Pinamonti, \emph{Loomis--Whitney inequalities in Heisenberg groups}, Math. Z. \textbf{301} (2022), 1983--2010.

\bibitem{FigalliMaggiPratelli}
A. Figalli, F. Maggi, and A. Pratelli, \emph{A mass transportation approach to quantitative isoperimetric inequalities}, Invent. Math. \textbf{182} (2010), 167--211.

\bibitem{FonsecaMuller}
I. Fonseca and S. M\"uller, \emph{A uniqueness proof for the Wulff Theorem}, Proc. Roy. Soc. Edinburgh Sect. A \textbf{119} (1991), 125--136.

\bibitem{FranceschiMontiRighiniSigalotti}
V. Franceschi, R. Monti, A. Righini, and M. Sigalotti, \emph{The isoperimetric problem for regular and crystalline norms in $\mathbb H^1$}, J. Geom. Anal. \textbf{33} (2023), Article 8.

\bibitem{Gromov}
M. Gromov, \emph{Groups of polynomial growth and expanding maps}, Inst. Hautes \'Etudes Sci. Publ. Math. \textbf{53} (1981), 53--78.

\bibitem{GromovEntropy}
M. Gromov, \emph{Entropy and isoperimetry for linear and non-linear group actions}, Groups Geom. Dyn. \textbf{2} (2008), 499--593.

\bibitem{MukherjeeExactMonotone}
M. Mukherjee, \emph{Exact versus monotone isoperimetric profiles of finitely generated groups}, preprint, \href{https://arxiv.org/abs/2509.22260v2}{arXiv:2509.22260v2}.

\bibitem{PansuPeriodic}
P. Pansu, \emph{Profil isop\'erim\'etrique, m\'etriques p\'eriodiques et formes d'\'equilibre des cristaux}, ESAIM Control Optim. Calc. Var. \textbf{4} (1999), 631--665.

\bibitem{PozueloExistence}
J. Pozuelo, \emph{Existence of isoperimetric regions in sub-Finsler nilpotent groups}, Anal. Geom. Metr. Spaces \textbf{13} (2025), Paper No.~20250025, 20 pp.

\bibitem{PozueloRitore}
J. Pozuelo and M. Ritor\'e, \emph{Pansu--Wulff shapes in $\mathbb H^1$}, Adv. Calc. Var. \textbf{16} (2023), 69--98.

\bibitem{RitoreCalibration}
M. Ritor\'e, \emph{A proof by calibration of an isoperimetric inequality in the Heisenberg group $\mathbb H^n$}, Calc. Var. Partial Differential Equations \textbf{44} (2012), 47--60.

\end{thebibliography}
\end{document}